\documentclass[review]{elsarticle}

\usepackage{theorem}
\usepackage{amssymb,amsmath}
\usepackage{bbm}
\usepackage{bm}

\usepackage{mathptmx}       % selects Times Roman as basic font
\usepackage{helvet}         % selects Helvetica as sans-serif font
\usepackage{courier}        % selects Courier as typewriter font
\usepackage{type1cm}        % activate if the above 3 fonts are
\usepackage{makeidx}         % allows index generation
\usepackage{graphicx}        % standard LaTeX graphics tool
\usepackage{multicol}        % used for the two-column index
\usepackage[bottom]{footmisc}% places footnotes at page bottom
\usepackage{subfigure}
\usepackage{booktabs}

\usepackage{epsfig} % for postscript graphics files
\usepackage{epstopdf}
\usepackage{color}
\usepackage{colordvi}

\newcommand{\norm}[2]         { \| {#1} \|_{#2} }                      % Norm
\newcommand{\bfm}[1]             { \mathbf{#1}     }             %
\newcommand{\SCZO}            { C^0(\Omega) }                    % C^0(Omega)
\newcommand{\SLTO}            { L^2(\Omega) }                       % L2(Omega)
\newcommand{\SLTOT}            { L^2(\Omega_T) }                       % L2(Omega)
\newcommand{\SHOO}            { H^1(\Omega) }                         % H1(Omega)
\newcommand{\SHOP}            { H^1(\Ph) }                              % broken H1 space
\newcommand{\VVK}            { {V(K_m)} }              % broken space of test functions
\newcommand{\VV}            { {V(\Ph)} }        % broken space of test functions
\newcommand{\UUUT}            { U(\Omega_T) }           %  space of solution u
\newcommand{\UUUhT}            { U^h(\Omega_T) }           %  space of solution u
\newcommand{\VVh}            { V^*(\Ph) }          % broken space of test functions
\newcommand{\SHOK}            { H^1(K_m) }                        % local H1 space
\newcommand{\SLTK}            { L^2(K_m) }                    % L2(K)
\newcommand{\SHMHGN}            { H^{-1/2}(\GN) }          % H^-1/2(Gamma_N)
\newcommand{\SHMHGD}            { H^{-1/2}(\GD) }          % H^-1/2(Gamma_N)
\newcommand{\SHHdK}            { H^{1/2}(\dKm) }         % H^1/2(boundary of Km)
\newcommand{\SHOOT}            { H^1(\Omega_T) }                   % H1(Omega)
\newcommand{\SHOOTT}            { H^1_{u_0}(\Omega_T) }                   % H1(Omega)
\newcommand{\SLOinf}            { L^\infty(\Omega)}          % scalar L infinity on Omega
\newcommand{\vn}             { \bfm{n}     }             % normal vector n
\newcommand{\GD}          { \Gamma_{D} }            % Boundary with Dirichlet Conditions
\newcommand{\GN}          { \Gamma_{N} }            % Boundary with Neumann Conditions
\newcommand{\Ph}            { \mathcal{P}_h }
\newcommand{\Kep}           { K_m \in \Ph}
\newcommand{\dKm}           { \partial K_m }

\newcommand{\dx}              { \; {\rm d} \bfm{x}   }                % dx (vector)
\newcommand{\dss}             { \, {\rm d} s   }                          % dx (vector)
\newcommand{\summa}[2]        { \overset{#2}{\underset{#1}{\sum}} } % summation_def for lower AND upper limits
\newcommand{\supp}[1]         { \underset{#1}{\sup} \, }        % supremum_def
\newcommand{\nn}              { \bfm{n} }                      % Vector n
\newcommand{\isdef}           { \overset{\text{def}}{=} } %  Equal sign used to define a quantity
\newcommand{\ds}              { \displaystyle }   %  Equal sign used to define a quantity

\theoremstyle{plain}
\theoremheaderfont{\bfseries \upshape}

\newtheorem{lem}{Lemma}[section]
\newtheorem{rem}{Remark}[section]

\begin{document}

\begin{frontmatter}
 \title{An Unconditionally Stable Space-Time FE Method for the Korteweg-de Vries Equation}

\author[1]{Eirik Valseth\corref{cor1}}
\ead{Eirik@utexas.edu}

\author[1]{Clint Dawson}
%\ead{Clint@oden.utexas.edu}

 \cortext[cor1]{Corresponding author}
% \fntext[fn1]{This is the first author footnote.}
% \fntext[fn2]{Another author footnote, this is a very long 
%   footnote and it should be a really long footnote. But this 
%   footnote is not yet sufficiently long enough to make two 
%   lines of footnote text.}

 \address[1]{Oden Institute for Computational Engineering and Sciences, The University of Texas at Austin, Austin, TX 78712, USA}


\biboptions{sort&compress}

%\newpageafter{title} 

%
%\maketitle
%
%

\begin{abstract}
We introduce an unconditionally stable finite element (FE) method, the automatic variationally stable FE (AVS-FE) method for the numerical analysis of the Korteweg-de Vries (KdV) equation. 
The AVS-FE method is a Petrov-Galerkin
method which employs the concept of optimal discontinuous test functions of
the discontinuous Petrov-Galerkin (DPG) method.
However, since AVS-FE method is a minimum residual method, we establish a global saddle point system instead of computing optimal test functions element-by-element. This system allows us to seek both the approximate solution 
of the KdV initial boundary value problem (IBVP) and a Riesz representer of the approximation
error. The AVS-FE method distinguishes itself from other minimum residual methods by using globally continuous
Hilbert spaces, such as $H^1$, while at the same time using broken Hilbert spaces for the test. 
Consequently, the AVS-FE approximations are classical $C^0$ continuous FE solutions.
The unconditional stability of this method allows us to solve the KdV equation space and time 
without having to satisfy a CFL condition. 
We present several numerical verifications for both linear and nonlinear versions of the KdV equation 
leading to optimal convergence behavior. Finally, we present a numerical verification of 
adaptive mesh refinements in both space and time for the nonlinear KdV equation.
\end{abstract}

\end{frontmatter}

\section{Introduction}
\label{sec:introduction}

The Korteweg-de Vries  equation~\cite{korteweg1895xli}, introduced in 1895 governs the propagation of one dimensional surface waves in water. There are also a multitude of interpretations
in which this equation governs a wide range of physical phenomena within wave propagation and these are summarized by Bhatta and Bhatti in~\cite{bhatta2006numerical}. 
This equation is a transient, third order nonlinear partial differential equation (PDE)
and has been 
analyzed by Holmer in~\cite{holmer2006initial}, where conditions for well-posedness of the differential operator are provided, i.e., conditions for its kernel being trivial. 
However, there are several issues requiring special consideration: $i)$ the nonlinearity of the KdV equation, 
$ii)$ the high order derivatives appearing in its differential operator, and $iii)$ its transient 
nature.
The latter being the most critical as FE methods for nonlinear problems have
been established in, e.g.,~\cite{oden2006finite} and mixed FE methods~\cite{BrezziMixed} can be employed to consider equivalent first order systems of PDEs.

The transient nature of the KdV leads to inherently unstable numerical approximations as first order 
(partial) time derivatives are advective transport terms. Thus, classical FE methods require 
mesh partitions that are sufficiently fine to achieve conditional stability and are therefore not 
suited for the KdV equation. Typically, transient PDEs are discretized using the method of lines 
to decouple spatial and temporal computations using FE methods in space and difference schemes 
in time. The numerical stability of the temporal discretization is then established by the Courant-Friedrichs-Lewy (CFL) condition~\cite{courant1928partiellen}. 
This type of approximation has been established for the KdV equation using multiple types of FE methods 
including discontinuous Galerkin (DG) method by multiple authors~\cite{levy2004local,yan2002local,yan2014existence}, the hybridized DG method of Samii \emph{et al.}~\cite{samii2016hybridized}, Galerkin methods using $C^1$ or higher bases~\cite{canivar2010taylor,amein2011small,bona1986fully,winther1980conservative,baker1983convergence}, and spectral element methods~\cite{llobell2020high,minjeaud2018high}. 
Due to the success of the aforementioned methods and the ease of implementation of a CFL condition for 
stability of time stepping schemes,  space-time FE methods have, to our best knowledge, not been 
applied to the KdV equation. In space-time FE methods, stabilized methods such as Streamlined Upwind 
Petrov-Galerkin by Hughes \emph{et al.}~\cite{Hughes1996,hughes1988space} can be applied to achieve conditional discrete 
stability. 
Increasing the dimension of the approximation space in FE methods 
increases the computational cost as the system of linear algebraic equations inevitably becomes larger. However, this cost is easily justified as space-time FE methods retain the attractive functional framework of FE methods 
in which \emph{a priori} and \emph{a posteriori} error estimates and adaptive strategies are 
available.

The AVS-FE method introduced by Calo, Romkes and Valseth in~\cite{CaloRomkesValseth2018} is an
 unconditionally stable Petrov-Galerkin FE method. As the 
 classical FE method, the AVS-FE method employs continuous trial spaces, whereas the test spaces are  
discontinuous. Thus, the AVS-FE method is a hybrid of the DPG 
method of Demkowicz and Gopalakrishnan~\cite{Demkowicz4, Demkowicz1, Demkowicz2, Demkowicz3, Demkowicz5, Demkowicz6}, as the 
test space consists of optimal discontinuous test functions, and the classical FE method. 
In addition to unconditional stability, the AVS-FE method satisfies a best approximation property and exhibits highly accurate flux predictions.

In this paper, we develop space-time AVS-FE approximations of the KdV equation. 
Following this introduction, we introduce the model KdV boundary value problem (BVP) and also notations and conventions in Section~\ref{sec:model_and_conv}.
Next, we review the AVS-FE methodology in Section~\ref{sec:AVS-review} and introduce the concepts
of Carstensen \emph{et al.}~\cite{carstensen2018nonlinear} to be employed to perform nonlinear iterations.
 In Section~\ref{sec:avs-fe}, we derive the equivalent AVS-FE weak 
formulation for the KdV BVP.
\emph{A priori} error estimates are introduced in Section~\ref{sec:avs-estimates}.
In Section~\ref{sec:verifications} we  perform multiple  numerical verifications for the KdV equation 
verifying numerical asymptotic convergence properties as well as a verification of an $h-$adaptive algorithm.
Finally, we conclude with remarks on the results and future works in Section~\ref{sec:conclusions}.

\section{KdV Equation and the AVS-FE Method }  
\label{sec:model_and_conv}

In this section, we present our model problem, i.e., the nonlinear KdV initial boundary value problem
(IBVP), present the notation and conventions we use, and give an overview of the AVS-FE method. 

\subsection{KdV Equation}
\label{sec:model_problem}
Let us consider the following KdV equation~\cite{korteweg1895xli}:
\begin{equation}  \label{eq:KDV_eq}
 \ds \frac{\partial u}{\partial t} + \frac{\partial}{\partial x}  (\beta u^2 + \alpha \frac{\partial^2 u}{\partial x^2}) = f(x,t),
\end{equation}
where $u$ represents the amplitude of a wave, the parameter $\beta$ is, as in~\cite{samii2016hybridized}, used to differentiate between the linear and nonlinear cases, 
and $\alpha$ denote the direction of the wave propagation.  
The PDE~\eqref{eq:KDV_eq} is posed on a domain $\Omega = \lbrack x_L, x_R \rbrack \subset \mathbb{R}^{1}$ and $t \in (0,T \rbrack$. Finally,
to establish a well defined KdV IBVP with a trivial kernel, the following set of initial and boundary conditions 
are considered as established by~\cite{holmer2006initial}, first we have the initial condition on $u$ and a boundary condition on the derivative of $u$:
\begin{equation} \label{eq:KDV_BCs}
\boxed{
\begin{array}{ll}
\ds u(x,0) & = u_0, \; x \in \Omega, \\
\ds \frac{\partial u}{\partial x} & = g_q, \; x = \, x_L \, \vee \, x_R, \\
 \end{array}}
\end{equation}
which are required. Next, one of the following sets of boundary conditions is needed to ensure the 
well prosedness of the IBVP:
\begin{equation} \label{eq:KDV_BCsII}
\boxed{
\begin{array}{ll}
\ds u(x,t) & = g_u, \; x = \, x_L \, \wedge \, x_R, \\
\ds u(x,t) & = g_u, \; x = \, x_L \wedge \frac{\partial^2 u(x,t)}{\partial x^2} \vn = g_p, \; x = \, x_R, \\
 \frac{\partial^2 u(x,t)}{\partial x^2} \vn & = g_p, \; x = \, x_L \wedge u(x,t) = g_u, \; x = \, x_L, \\
 \end{array}}
\end{equation}
where $\vn$ denote the unit scalar, i.e., $\vn = \pm 1$.

We restrict ourselves to the first case in~\eqref{eq:KDV_BCsII} for the sake of brevity, and label the boundary where we have conditions on $u$ by $\Gamma_D$ and the boundary 
where we have a condition on the derivative $\frac{\partial u}{\partial x}$ by $\Gamma_N$.
Thus, our  KdV model IBVP is:
\begin{equation} \label{eq:KDV_IBVP}
\boxed{
\begin{array}{l}
\text{Find }  u  \text{ such that:}    
\\[0.05in] 
\qquad 
\begin{array}{rrl}
 \ds \frac{\partial u}{\partial t} + \frac{\partial}{\partial x}  (\beta u^2 + \alpha \frac{\partial^2 u}{\partial x^2}) & = &  f,  \text{ in } \, \Omega\times(0,T], 
 \\[0.05in]
 u & = & u_0 \text{ on } \Omega \times (t = 0),   \\ 
u & = & g_{u}, \text{ on } \Gamma_D,  \\  
\frac{\partial u}{\partial x} & = & g_{q}, \text{ on } \Gamma.  \\  
 \end{array}
 \end{array}
}
\end{equation}

\subsection{Review of the AVS-FE Method}
\label{sec:AVS-review}
The AVS-FE method~\cite{CaloRomkesValseth2018,valseth2020CH,eirik2019thesis,valseth2020goal} is a 
conforming and unconditionally stable FE method that employs the concept of optimal test spaces.
This method remains attractive in particular due to its unconditional stability property, regardless 
of the differential operator, as well as highly accurate flux, or derivative, approximations. 
The method is a Petrov-Galerkin method in which the trial space consists of classical,
globally continuous, Hilbert spaces, whereas the test space consists of broken Hilbert spaces. Thus, 
the test functions are square integrable functions globally that are allowed higher order regularity on 
each element in the FE mesh. In particular, this test space is spanned by a basis that is optimal, 
in the sense that it leads to unconditional stability by computing on-the-fly optimal test functions 
in the spirit of the discontinuous Petrov-Galerkin (DPG) method~\cite{Demkowicz4, Demkowicz1, Demkowicz2, Demkowicz3, Demkowicz5, Demkowicz6}. In the following review, we present key features of the 
AVS-FE method and omit some details, a thorough introduction can be found in ~\cite{CaloRomkesValseth2018} for the AVS-FE method and~\cite{ValsethBehnoud2020a} for the space-time version of the AVS-FE method.

To introduce the AVS-FE method we consider the abstract weak formulation:
\begin{equation} \label{eq:abstract_weak_form_AVS}
\boxed{
\begin{array}{ll}
\text{Find } \mathbbm{u} \in U & \hspace{-0.05in} \text{ such that:}
\\[0.05in]
 &  \quad B(\mathbbm{u};\mathbbm{v}) = F(\mathbbm{v}), \quad \forall \mathbbm{v} \in V, 
 \end{array}}
\end{equation}
where $\mathbbm{u}$ and $\mathbbm{v}$ are the trial and test functions, respectively, $U$ is the trial space, $V$ the (broken) test
space, $B:U\times V\longrightarrow \mathbb{R}$ is the bilinear (or sesquilinear) form, and 
  $F:V \longrightarrow \mathbb{R}$  the linear 'load' functional.
%  \, and $\Ph$ denotes the partition of $\Omega$ into elements (see~\eqref{eq:domain}).
The kernel of the underlying boundary value problem (BVP) is assumed to be trivial to ensure the 
uniqueness of the weak solutions of~\eqref{eq:abstract_weak_form_AVS} and $B(\cdot,\cdot)$ and $F(\cdot)$
satisfy continuity conditions. 
Thus, the remaining item for well-posedness of~\eqref{eq:abstract_weak_form_AVS} is the \emph{inf-sup}
condition:
\begin{equation}
\label{eq:abstract_inf_sup}
 \gamma > 0 :
 \supp{\mathbbm{v}\in V\setminus \{0\}} 
     \frac{|B(\mathbbm{u};\mathbbm{v})|}{\norm{\mathbbm{v}}{V}} \ge \gamma \norm{\mathbbm{u}}{U}.
\end{equation}
Generally, satisfying this condition becomes problematic in the discrete case, where the supremum over 
the test space $V$ is not identical to the supremum over a (discrete) subspace $V^h \subset V$.
The AVS-FE method ensures satisfaction of both continuous and discrete \emph{inf-sup} conditions by 
employing a test space that ensures $\gamma = \gamma^h = 1$, as well as an 
alternative norm on the trial space.
The kernel of the underlying BVP being trivial lets
 us introduce the following \emph{energy norm} $\norm{\cdot}{B}: U\longrightarrow [0,\infty)$:
\begin{equation}
\label{eq:abstract_energy_norm}
\norm{\mathbbm{u}}{B} \isdef \supp{\mathbbm{v}\in V\setminus \{\mathbf{0}\}} 
     \frac{|B(\mathbbm{u};\mathbbm{v})|}{\norm{\mathbbm{v}}{V}}.
\end{equation}

Following the philosophy of the DPG method, the optimal test space is spanned by functions that are
solutions of a global Riesz representation problem:
\begin{equation} \label{eq:abstract_riesz_problem}
\begin{array}{rcll}
\ds \left(\, \mathbbm{p};\mathbbm{v} \, \right)_V &  \! \! =  \! & B(\,\mathbbm{u};\mathbbm{v} \, ),& \, \forall \mathbbm{v}\in V, 
\end{array}
\end{equation}
for each $u \in U$, where the operator on the left hand side (LHS) is an inner product. 
Clearly, the Riesz problem is well posed and its solutions are to be ascertained numerically by a FE approximation. Thankfully, as the test space is broken, the
approximations of~\eqref{eq:abstract_riesz_problem} are computed in a decoupled fashion, element-by-element. We note that the solutions of~\eqref{eq:abstract_riesz_problem} may exhibit 
boundary layers on each element depending on the inner product in the LHS and the action 
of the bilinear form onto functions in the trial space. However, in the AVS-FE method, the inner product is 
defined solely by the regularity of the test space and the solutions of~\eqref{eq:abstract_riesz_problem}  do not exhibit boundary layers~\cite{eirik2019thesis}.
Due to the Riesz representation problem, the continuous and discrete
well posedness of~\eqref{eq:abstract_weak_form_AVS} follows from the Generalized Lax-Milgram Theorem, as $B(\cdot,\cdot)$ 
satisfies both the inf-sup condition as well as the continuity condition in terms of the energy
norm~\eqref{eq:abstract_energy_norm} (see~\cite{eirik2019thesis} or~\cite{Demkowicz5,Demkowicz4} for details).

The derivation of weak formulations for the AVS-FE method, largely follows the approach of classical mixed FE methods in that the trial space consists of  global Hilbert spaces. However, it differs significantly 
from these methods as the test space is a broken Hilbert space as in discontinuous Galerkin or 
DPG methods. Consequently, in the FE discretization of the AVS-FE weak formulation, we employ 
classical \emph{continuous} FE bases such as Lagrange or Raviart-Thomas for the trial space, whereas 
the test space is spanned by the \emph{discontinuous} optimal test functions computed from 
the discrete analogue of~\eqref{eq:abstract_riesz_problem}. Hence,
the FE discretization of~\eqref{eq:abstract_weak_form_AVS}  governing the approximation $\mathbbm{u}^h\in U^h$ of $\mathbbm{u}$ is:
\begin{equation} \label{eq:discrete_form}
\boxed{
\begin{array}{ll}
\text{Find} &  \mathbbm{u}^h \in U^h \; \text{ such that:}
\\[0.1in]
 &   B(\mathbbm{u}^h;\mathbbm{v}^h) = F(\mathbbm{v}^h), \quad \forall \mathbbm{v}^h\in V^{*}, 
 \end{array}}
\end{equation}
where the finite dimensional subspace of test functions $V^{*} \subset V$ is spanned by the optimal test functions.

The DPG philosophy used to construct $V^{*}$ ensures that the
discrete problem~\eqref{eq:discrete_form} inherits the continuity and inf-sup constants
of the continuous problem. 
Hence, the AVS-FE discretization is unconditionally stable for any choice of element size $h$ and local degree of polynomial approximation $p$. Furthermore, the 
global stiffness matrix is symmetric and positive definite regardless of the character of the underlying differential operator.
The globally continuous trial space also has the consequence that the optimal test functions have identical support to its corresponding trial function, i.e., compact. Consequently, the bandwidth of the 
global stiffens matrix is the same as in classical mixed FE methods. 
Clearly, the cost of assembling the global system of linear algebraic equations is greater than in 
classical FE methods due to the solution of the local Riesz representation problems. However, 
this cost is kept to a minimum as these local problems are computed at the same degree of approximation 
as the trial functions and the cost of computing the optimal test functions is incurred at the element level as part of the FE assembly process. 

\begin{rem} \label{rem:mixed}
As an alternative to computing optimal test functions on-the-fly to construct
the FE system of linear algebraic equations, one can 
consider another equivalent, interpretation of the AVS-FE method. 
This alternative interpretation is in the DPG
community often referred to as a mixed or saddle point problem and
results from a constrained minimization interpretation of the DPG and AVS-FE methods 
in which the Riesz representation problem is used to define a constraint equation to the weak form:
\begin{equation} \label{eq:constrained_min_problem}
\boxed{
\begin{array}{rl}
\text{Find } \mathbbm{u}^h \in  U^h, \mathbbm{\hat{e}}^h \in V^h  & \hspace{-0.15in} \text{ such that:}
\\[0.05in]
   \quad \left(\, \mathbbm{\hat{e}}^h,\mathbbm{v}^h \, \right)_{V} - B(\mathbbm{u}^h,\mathbbm{v}^h) & =  - F(\mathbbm{v}^h), \quad \forall \mathbbm{v}^h \in V^h,  \\
  \quad B'_{\mathbbm{u}}(\mathbbm{p}^h,\mathbbm{\hat{e}}^h)& =  0, \quad \forall \, \mathbbm{p}^h \in U^h.
 \end{array}}
\end{equation}
Where the second equation represents a 
constraint in which the Gateaux derivative of the bilinear form is acting on the approximate "error representation" function $\mathbbm{\hat{e}}^h$. This function is 
a Riesz representer of the approximation error $\mathbbm{u}-\mathbbm{u}^h$. 
The energy norm of this approximation error is identical to the norm of the error representation 
function due to~\eqref{eq:abstract_riesz_problem}.
Thus, the 
norm of the approximate error representation function $\norm{\mathbbm{\hat{e}}^h}{V}$ is an a posteriori error estimate. For details on these error indicators and the derivation of the mixed formulation, see~\cite{eirik2019thesis} or~\cite{demkowicz2014overview}.
Note for differential 
operators which are linear, the Gateaux derivative of the bilinear form is identical 
to itself. Hence, for nonlinear differential operators such as the KdV equation, the Gateaux derivative $B'_{\mathbbm{u}}(\cdot,\cdot)$ of $B(\cdot,\cdot)$ must be established. 

The  cost of solving the resulting system of linear algebraic equations from~\eqref{eq:constrained_min_problem} is larger than the 'classical' AVS-FE method since now the optimal test functions are essentially computed by solving global problems. However, it has a clear advantage for mesh adaptive strategies, since upon 
solving~\eqref{eq:constrained_min_problem}, it provides a posteriori error 
estimators and error indicators that can drive the mesh adaptive process. Furthermore the efforts 
required in the implementation of the method are small into high level FE solvers such FEniCS~\cite{alnaes2015fenics}.
\end{rem}

\section{AVS-FE Weak Formulation and Discretization of The KdV Equation}  
\label{sec:avs-fe}

With the notations introduced in Section~\ref{sec:model_and_conv} and the review of the
AVS-FE method above, we proceed to derive the AVS-FE weak formulation for the KdV IBVP~\eqref{eq:KDV_IBVP}:
\begin{equation*} 
\boxed{
\begin{array}{l}
\text{Find }  u  \text{ such that:}    
\\[0.05in] 
\qquad 
\begin{array}{rrl}
 \ds \frac{\partial u}{\partial t} + \frac{\partial}{\partial x}  (\beta u^2 + \alpha \frac{\partial^2 u}{\partial x^2}) & = &  f,  \text{ in } \, \Omega_T, 
 \\[0.05in]
 u & = & u_0 \text{ on } \Omega \times (t = 0),   \\ 
u & = & g_{u}, \text{ on } \Gamma_D,  \\  
\frac{\partial u}{\partial x} & = & g_{q}, \text{ on } \Gamma_N,  \\  
 \end{array}
 \end{array}
}
\end{equation*}
where $\Omega_T = \Omega\times(0,T]$ is the space-time domain and both $\beta$ and $\alpha$ belong to
 $\SLOinf$. We assume that the source $f \in \SLTO$, but this assumption is not strictly necessary for 
 well-posedness of the AVS-FE weak formulation.
The starting point of the derivation is the regular partition 
$\Ph$ of $\Omega_T$ into elements $K_m$, such that:
\begin{equation}
\label{eq:domain}
  \Omega_T = \text{int} ( \bigcup_{\Kep} \overline{K_m} ).
\end{equation}
We apply a mixed FE methodology and introduce two auxiliary variables:
\begin{itemize}
\item $q = \frac{\partial u}{\partial x}$.
\item $p = \frac{\partial q}{\partial x}$.
\end{itemize}
The first-order system of the KdV IBVP~\eqref{eq:KDV_IBVP} is therefore: 
\begin{equation} \label{eq:KDV_IBVP_first_order}
\boxed{
\begin{array}{l}
\text{Find }  (u,q,p)  \text{ such that:}    
\\[0.05in] 
\qquad 
\begin{array}{rcl}
\ds   q - \frac{\partial u}{\partial x} & =  & 0, \quad \text{ in } \, \Omega, 
  \\[0.055in]
  \ds   p - \frac{\partial q}{\partial x} & =  & 0, \quad \text{ in } \, \Omega, 
  \\[0.045in]
  \ds  \frac{\partial u}{\partial t} + \beta \frac{\partial (\beta u^2) }{\partial x} + \alpha \frac{\partial p}{\partial x} & =  & f, \quad \text{ in } \, \Omega_T, 
 \\[0.025in]
 \qquad u & = & u_0, \quad \text{ on } \Omega \times (t = 0), \\
 \qquad u & = & g_{u}, \quad \text{ on } \Gamma_D,  \\  
 \qquad q & = & g_{q}, \quad \text{ on } \Gamma_N,  \\  
 \end{array}
 \end{array}
}
\end{equation}

Next, enforce the PDEs~\eqref{eq:KDV_IBVP_first_order} weakly on each
element $\Kep$, i.e.,
\begin{equation} \label{eq:weak_IBVP_L2_KDV}
\begin{array}{c}
\text{Find } \; (u,q,p) \in \SHOOT\times\SHOOT\times\SHOOT:  
\\[0.1in] 
\ds   \int_{K_m} \biggl\{\left[ q - \frac{\partial u}{\partial x}\right]  s_m \, 
\ds+ \left[ p - \frac{\partial q}{\partial x} \right] w_m \,
\ds+ \left[ \frac{\partial u}{\partial t} + \beta \frac{\partial (\beta u^2) }{\partial x} + \alpha \frac{\partial p}{\partial x}  \right] \, v_m  \biggr\} \dx   =  
\ds   \int_{K_m} f\, v_m \dx, 
 \\ \hspace{3.5in}  \forall (v_m,w_m,s_m) \in \SLTK\times\SLTK\times\SLTK,
  \end{array}
\end{equation}
where $\dx = \text{d}x \, \text{d}t$.
We proceed to apply Green's identity to all terms with spatial partial derivatives:
\begin{equation} \label{eq:ultraweak_IBVP_KDV}
\begin{array}{c}
\text{Find } \; (u,q,p) \in \SHOOT\times\SHOOT\times\SHOOT:  
\\[0.1in] 
\ds  \int_{K_m} \biggl\{ q \, s_m + \frac{\partial s_m}{\partial x} u \, 
\ds+  p\, w_m + \frac{\partial w_m}{\partial x}   q_m \,
\ds+ \frac{\partial u}{\partial t} \, v_m - \beta \frac{ \partial (v_m) }{\partial x} \, u^2 - \alpha \frac{\partial v_m}{\partial x} \, p \, \biggr\} \dx  \\
\ds + \oint_{\dKm} \biggl\{ \alpha \,  \gamma^m_\nn(p) \, \gamma^m_0(v_m) + \beta \,  \gamma^m_\nn(u^2) \, \gamma^m_0(v_m) \, - \,  \gamma^m_\nn(q) \, \gamma^m_0(w_m) \,- \,  \gamma^m_\nn(u) \, \gamma^m_0(s_m)  \, \biggr\} \dss
  =  
\ds \int_{K_m} f\, v_m \dx, 
 \\ \hspace{3.1in} \forall (v_m,w_m,s_m) \in \SHOK\times\SHOK\times\SHOK.
  \end{array}
\end{equation}
Note that, since we are in one dimension, the boundary integrals are simply function evaluations 
on the left or right hand side of the elements.
We also employ engineering convention here and use an integral representation of the boundary terms 
which are to be interpreted as duality pairings on $\dKm$,
the operators $\gamma^m_0: \SHOK: \longrightarrow \SHHdK$ and $\gamma^m_\nn: \SHOK: \longrightarrow \SHHdK$, i.e., $\gamma^m_\nn(\cdot) = \gamma^m_0(\cdot \, \nn)$, denote trace operators (e.g., see~\cite{Girault1986}) on $K_m$; and $\nn$ is the unit normal scalar on  the element boundary $\dKm$ of $K_m$.

In an effort to enforce all boundary conditions in a weak manner, we decompose the boundary terms 
including $u$ and $q$ into portions intersecting the global boundaries $\Gamma_D$ and $\Gamma_N$:
\begin{equation} \label{eq:ultraweak_IBVP_KDV_decomp}
\begin{array}{c}
\text{Find } \; (u,q,p) \in \SHOOT\times\SHOOT\times\SHOOT:  
\\[0.1in] 
\ds \int_{K_m} \biggl\{ q \, s_m + \frac{\partial s_m}{\partial x} u \, 
\ds+  p\, w_m + \frac{\partial w_m}{\partial x}   q_m \,
\ds+ \frac{\partial u}{\partial t} \, v_m - \beta \frac{ \partial (v_m) }{\partial x} \, u^2 - \alpha \frac{\partial v_m}{\partial x} \, p \, \biggr\} \dx  \\
\ds + \oint_{\dKm} \biggl\{ \alpha \,  \gamma^m_\nn(p) \, \gamma^m_0(v_m) + \beta \,  \gamma^m_\nn(u^2) \, \gamma^m_0(v_m)   \, \biggr\} \dss   \ds - \oint_{\dKm\setminus \GD } \biggl\{    \gamma^m_\nn(u) \, \gamma^m_0(s_m)  \, \biggr\} \dss  \\ \ds -  \oint_{\dKm\setminus \GN } \biggl\{   \gamma^m_\nn(q) \, \gamma^m_0(w_m) \, \biggr\} \dss 
\ds - \oint_{\dKm \cap \GD } \biggl\{  \gamma^m_\nn(u) \, \gamma^m_0(s_m)  \, \biggr\} \dss
- \oint_{\dKm \cap \GN } \biggl\{  \gamma^m_\nn(q) \, \gamma^m_0(w_m)  \, \biggr\} \dss \\[0.1in]
\hspace{1in}  =  \ds  \int_{K_m} f\, v_m \dx,  \quad  \forall (v_m,w_m,s_m) \in \SHOK\times\SHOK\times\SHOK.
  \end{array}
\end{equation}
The boundary conditions are subsequently enforced weakly on $\GD$ and $\GN$, the initial condition 
is enforced in a strong manner and is incorporated into the trial space. Then,
 we sum contributions
 from all local integral statements $\Kep$, and arrive at the global weak formulation:
\begin{equation} \label{eq:ultraweak_IBVP_KDV_decomp2}
\begin{array}{c}
\text{Find } \; (u,q,p) \in \UUUT:  
\\[0.1in] 
\ds \summa{\Kep}{} \int_{K_m} \biggl\{ q \, s_m + \frac{\partial s_m}{\partial x} u \, 
\ds+  p\, w_m + \frac{\partial w_m}{\partial x}   q_m \,
\ds+ \frac{\partial u}{\partial t} \, v_m - \beta \frac{ \partial (v_m) }{\partial x} \, u^2 - \alpha \frac{\partial v_m}{\partial x} \, p \, \biggr\} \dx  \\
\ds + \oint_{\dKm} \biggl\{ \alpha \,  \gamma^m_\nn(p) \, \gamma^m_0(v_m) + \beta \,  \gamma^m_\nn(u^2) \, \gamma^m_0(v_m)   \, \biggr\} \dss   \ds - \oint_{\dKm\setminus \GD } \biggl\{    \gamma^m_\nn(u) \, \gamma^m_0(s_m)  \, \biggr\} \dss  \\
\ds - \oint_{\dKm \setminus \GN } \biggl\{  \gamma^m_\nn(q) \, \gamma^m_0(w_m)  \, \biggr\} \dss \\[0.1in]
  =  \ds \summa{\Kep}{} \int_{K_m} f\, v_m \dx + \oint_{\dKm \cap \GD } \biggl\{  g_u \, \gamma^m_0(s_m)  \, \biggr\} \dss +  \oint_{\dKm\cap \GN } \biggl\{   g_q \, \gamma^m_0(w_m) \, \biggr\} \dss ,  \\ \hspace{3.5in}  \forall (v_m,w_m,s_m) \in \VV,
  \end{array}
\end{equation}
where the trial and test spaces are defined:
\begin{equation}
\label{eq:function_spaces_KDV}
\begin{array}{c}
\ds \UUUT \isdef \SHOOTT\times\SHOOT\times\SHOOT ,
\\[0.15in]
\ds \VV \isdef  \SHOP\times\SHOP\times\SHOP,
\end{array}
\end{equation}
with $\SHOOT$ being classical first order Hilbert spaces over $\Omega_T$ and:
\begin{equation}
\label{eq:SHOOTT_SPACE}
\SHOOTT \isdef \biggl\{ v\in\SHOOT: \quad v|_{\Omega \times (t = 0)} = u_0 \biggr\}.
\end{equation}
\begin{equation}
\label{eq:broken_h1_space}
\SHOP \isdef \biggl\{ v\in\SLTO \times (0,T]: \quad v_m \in \SHOK, \; \forall \Kep\biggr\}.
\end{equation}
We also define the norms $\norm{\cdot}{\UUUT}:  \UUUT \!\! \longrightarrow \!\! [0,\infty)$ and $\norm{\cdot}{\VV}: \VV\! \! \longrightarrow\! \! [0,\infty)$ as:
\begin{equation}
\label{eq:broken_norms_KDV}
\begin{array}{l}
\ds \norm{(u,q,p}{\UUUT} \isdef \sqrt{ \biggl[ (u,u)_{\SHOOT} + (q,q)_{\SHOO}  +  (p,p)_{\SHOO}  \biggr] },
\\[0.2in]
\ds   \norm{(v,w,s)}{\VV} \isdef  \ds \sqrt{\summa{\Kep}{}\int_{K_m} \biggl[  h_m^2  \frac{\partial v_m}{\partial x}^2  + v_m^2 + h_m^2 \frac{\partial w_m}{\partial x}^2 + w_m^2    +  h_m^2 \frac{\partial s_m}{\partial x}^2 + s_m^2\biggr] \dx },
 \end{array}
\end{equation}
where $h_m$ denotes the element diameter.

\begin{rem}
The definition of the inner product that induces the norm on $\VV$ in~\eqref{eq:broken_norms_KDV}, defines the inner product in the LHS of the Riesz representation problem
governing the optimal test functions (see~\eqref{eq:abstract_riesz_problem}). 
This inner product follows naturally from our derivation of the weak formulation as it is simply the 
broken norm on the test space. The factors $h_m^2$ are needed to ensure that as $h_m \rightarrow 0$, the norm $\norm{\cdot}{\VV}$ remains bounded. 
Consequently, the resulting optimal test functions can be thought of as solutions to a balanced 
reaction-diffusion problem thereby ensuring that the optimal test functions do not exhibit 
local boundary layers.
\end{rem}

By introducing the sesquilinear and linear forms $B:\UUUT\times\VV\longrightarrow \mathbb{R}$ and $F:\VV\longrightarrow \mathbb{R}$:
\begin{equation} \label{eq:B_and_F_KDV}
\begin{array}{c}
\hspace{-0.15in} B((u,q,p);(v,w,s)) \isdef 
\ds \summa{\Kep}{} \biggl\{ \int_{K_m} \biggl\{ q \, s_m + \frac{\partial s_m}{\partial x} u \, 
\ds+  p\, w_m + \frac{\partial w_m}{\partial x}   q \,
\ds+ \frac{\partial u}{\partial t} \, v_m - \beta \frac{ \partial (v_m) }{\partial x} \, u^2 - \alpha \frac{\partial v_m}{\partial x} \, p \, \biggr\} \dx  \\
\ds + \oint_{\dKm} \biggl\{ \alpha \,  \gamma^m_\nn(p) \, \gamma^m_0(v_m) + \beta \,  \gamma^m_\nn(u^2) \, \gamma^m_0(v_m)   \, \biggr\} \dss   \ds - \oint_{\dKm\setminus \GD } \biggl\{    \gamma^m_\nn(u) \, \gamma^m_0(s_m)  \, \biggr\} \dss  \\
\ds - \oint_{\dKm \setminus \GN } \biggl\{  \gamma^m_\nn(q) \, \gamma^m_0(w_m)  \, \biggr\} \dss\biggr\}, \\ 
 F(v,w,s) \isdef \ds \summa{\Kep}{} \biggl\{ \int_{K_m} f\, v_m \dx + \oint_{\dKm \cap \GD } \biggl\{  g_u \, \gamma^m_0(s_m)  \, \biggr\} \dss +  \oint_{\dKm \cap \GN } \biggl\{   g_q \, \gamma^m_0(w_m) \, \biggr\} \dss   \biggr\},
 \end{array}
\end{equation}
we can write the weak formulation~\eqref{eq:ultraweak_IBVP_KDV_decomp2} compactly:
\begin{equation} \label{eq:weak_form_KDV}
\boxed{
\begin{array}{ll}
\text{Find } (u,q,p) \in \UUUT   \text{ such that:}
\\[0.05in]
   \quad B((u,q,p);(v,w,s)) = F(v,w,s), \quad \forall (v,w,s) \in \VV. 
 \end{array}}
\end{equation}
Now, we can introduce the corresponding energy norm and Riesz representation problems
(see \eqref{eq:abstract_energy_norm}) for the 
AVS-FE weak form of the KdV IBVP~\eqref{eq:weak_form_KDV}. Hence, 
we have the following well-posedness result:
\begin{lem}
\label{lem:well_posed_cont}
Let $f \in (\SHOP)'$, and the boundary data $g_u \in \SHMHGD$ and $g_q \in \SHMHGN$. Then, the weak formulation~\eqref{eq:weak_form_KDV} has a unique solution and is well posed.
\end{lem} 
\emph{Proof}: The kernel of the KdV IBVP is trivial for the 
chosen boundary conditions. Consequently, the kernel of the sesquilinear form in~\eqref{eq:weak_form_KDV} is also trivial. Then, the Generalized Lax-Milgram Theorem~\cite{babuvska197finite} ensures that the weak statement~\eqref{eq:weak_form_KDV} is 
well posed in terms of the energy norm~\eqref{eq:abstract_energy_norm} with inf-sup and continuity constants equal to unity.
\newline \noindent ~\qed

\noindent The AVS-FE weak formulation~\eqref{eq:weak_form_KDV} essentially represents a DPG ultraweak formulation
(see, e.g.,~\cite{Demkowicz5} in which only the test space is broken. 
Note that this weak formulation is not a unique choice and multiple other choices of where to apply  Green's Identity are possible.

\subsection{AVS-FE Discretizations}  
\label{sec:avs-discretization}

As in classical FE methods,
to seek numerical approximations $(u^h,q^h,p^h)$ of $(u,q,p)$, the AVS-FE method employs a discrete trial space $\UUUhT \subset \UUUT$ that consists of classical FE basis functions. This is a major advantage of the AVS-FE method, as the derivation of the 
weak formulation is performed such that the corresponding FE discretization is as close to the classical FE method as possible. Here, since the trial space consist of three $H^1$ Hilbert spaces, 
we employ standard $\SCZO$ continuous Lagrange polynomials for all trial variables. 
However, to ensure the discrete stability of the space-time discretization, we use the philosophy of
the DPG method to construct the test space.
Hence, the AVS-FE discretization of~\eqref{eq:weak_form_KDV} is: 
\begin{equation} \label{eq:discrete_form_KDV}
\boxed{
\begin{array}{ll}
\text{Find } (u^h,q^h,p^h) \in \UUUhT   \text{ such that:}
\\[0.05in]
   \quad B((u^h,q^h,p^h);(v^h,w^h,s^h)) = F(v^h,w^h,s^h), \quad \forall (v^h,w^h,s^h) \in \VVh, 
 \end{array}}
\end{equation}
where the finite dimensional test space $\VVh \subset \VV$ is spanned by numerical
approximations of the test optimal functions through the Riesz representation problems. 
The optimal test functions ensure the well posedness of~\eqref{eq:discrete_form_KDV} 
and we can employ the same arguments of Lemma~\ref{lem:well_posed_cont} to show this.
The dicretization~\eqref{eq:discrete_form_KDV} governs the space-time solutions $(u^h,q^h,p^h)$ and 
since it remains well posed for any choice of mesh parameters we do not require a CFL condition 
to ensure stability. 

While the AVS-FE method guarantees the stability of~\eqref{eq:discrete_form_KDV}, we have not
addressed its nonlinearity. 
Normally, one would employ a linearization technique and an iterative procedure to establish solutions of~\eqref{eq:discrete_form_KDV}. However, we shall take the approach introduced in Remark~\ref{rem:mixed}, and establish a saddle point problem in which the linearization  occurs naturally
in the form of a constraint (see~\cite{carstensen2018nonlinear} for details).
Thus, we consider the equivalent saddle point problem 
 in which we seek both $(u^h,q^h,p^h)$ 
 and the function $(\psi^h, \varphi^h, \xi^h)$, which is called the  error representation function: 
\begin{equation} \label{eq:disc_mixed}
\boxed{
\begin{array}{ll}
\text{Find } (u^h,q^h,p^h) \in  \UUUhT, (\psi^h, \varphi^h, \xi^h) \in V^h(\Ph)  \text{ such that:}
\\[0.05in]
   \quad \left(\, (\psi^h, \varphi^h, \xi^h),(v^h,w^h,s^h) \, \right)_\VV - B((u^h,q^h,p^h);(v^h,w^h, s^h)) & =  -F(v^h,w^h,s^h), \\& \quad  \hspace*{-1in}\forall (v^h,w^h,s^h) \in V^h(\Ph) ,  \\
  \quad B'_{\mathbbm{u}}((a^h,b^h,c^h));(\psi^h, \varphi^h, \xi^h) & =  0, \\& \quad \hspace*{-1in}\forall \, (a^h,b^h,c^h) \in \UUUhT.
 \end{array}}
\end{equation}
Where the operator $B'_{\mathbbm{u}}:\UUUT\times\VV\longrightarrow \mathbb{R}$ is the 
first order Gateaux derivative of the sesquilinear form $B$ with respect to  
$\mathbbm{u} = (u,q,p)$::
\begin{equation} \label{eq:B_prime_KDV}
\begin{array}{c}
B'_{\mathbbm{u}}((a,b,c);(v, w, s )) \isdef 
\ds \summa{\Kep}{} \biggl\{ \int_{K_m} \biggl\{ b \, s_m + \frac{\partial s_m}{\partial x} a \, 
\ds+  c\, w_m + \frac{\partial w_m}{\partial x}   b \,
\ds+ \frac{\partial a}{\partial t} \, v_m - 2 \beta \frac{ \partial (v_m) }{\partial x} \, ua - \alpha \frac{\partial v_m}{\partial x} \, c \, \biggr\} \dx  \\
\ds + \oint_{\dKm} \biggl\{ \alpha \,  \gamma^m_\nn(c) \, \gamma^m_0(v_m) + 2 \beta \,  \gamma^m_\nn(ua) \, \gamma^m_0(v_m)   \, \biggr\} \dss   \ds - \oint_{\dKm\setminus \GD } \biggl\{    \gamma^m_\nn(a) \, \gamma^m_0(s_m)  \, \biggr\} \dss  \\
\ds - \oint_{\dKm \setminus \GN } \biggl\{  \gamma^m_\nn(b) \, \gamma^m_0(w_m)  \, \biggr\} \dss\biggr\},
 \\[0.15in]
 \end{array}
\end{equation}

Essentially, the error representation function is the solution of a Riesz representation problem (see~\eqref{eq:abstract_riesz_problem}) where the RHS is the residual functional from the FE 
approximation of the trial functions:
\begin{equation} \label{eq:Ries_KDV_Error}
\begin{array}{rcll}
\ds \left(\, (\psi, \varphi, \xi);(v,w,s) \, \right)_{\VV} &  \! \! =  \! & B(\,(u^h,q^h,p^h);(v,w,s) \, ) - F(v,w,s), \quad \forall (v,w,s) \in \VV.
\end{array}
\end{equation}
Hence, the error representation function is an exact representation of the approximation error in terms 
of the energy norm:
\begin{equation}
\label{eq:energynorm_eqv}
\norm{(u-u^h,q-q^h,p-p^h)}{B} = \norm{(\psi, \varphi, \xi)}{\VV},
\end{equation}
and its approximation is an \emph{a posteriori} error estimate (see~\cite{Demkowicz2} for a
thorough introduction). Furthermore, since the space $\VV$ is broken, the computation of the 
norm can be performed element-wise and its local restriction is an error indicator:
\begin{equation}
\label{eq:enrr_ind}
\eta = \norm{(\psi^h, \varphi^h, \xi^h)}{\VVK}.
\end{equation}
This error indicator has been successfully applied to a wide range of problems in both the DPG and AVS-FE methods~\cite{carstensen2016breaking,roberts2015discontinuous,valseth2020CH}.

\subsection{A Priori Error Estimates}  
\label{sec:avs-estimates}

In this section, we present \emph{a priori} error estimates for the AVS-FE method applied to the 
linear version of the 
KdV equation in terms of norms of the approximation error. Extensive proofs are not provided here for the sake of brevity but rather outlines 
and references to the required literature are given where applicable. 
Here, the norms we are interested in are the energy norm, and the Sobolev norms on $\SHOOT$ and $\SLTOT$.  

First, we present the following lemma of the \emph{a priori} error estimate in terms of the energy norm:
\begin{lem} \label{lem:energy_error}
Let $(u,q,p) \in \UUUT$ be the exact solution of the AVS-FE weak formulation~\eqref{eq:weak_form_KDV} and  $(u^h,q^h,p^h) \in \UUUhT$ its corresponding AVS-FE approximation. Then:
\begin{equation} \label{eq:energy_rate}
\ds \exists \, C > 0 \, : \norm{(u-u^h,q-q^h,p-p^h)}{B} \leq  C\, \ds h^{\,\mu-1},
\end{equation} 
where $h$ is the maximum element diameter, $\mu = \text{ min } (p_u+1,r)$,  $p_u$ the minimum polynomial degree of approximation of $u^h$
in the mesh, and  $r$ the regularity of the solution $u$ of the governing PDE~\eqref{eq:KDV_eq}. 
\end{lem} 
\emph{Proof:}
The bound~\eqref{eq:energy_rate} is a consequence of:
\begin{itemize}
\item The  best approximation property of the AVS-FE method in terms of the 
energy norm~\eqref{eq:abstract_energy_norm} (see Lemma 2.3.2 in~\cite{eirik2019thesis});
\item The norm equivalence in~\eqref{eq:energynorm_eqv};
\item The existence of polynomial interpolation operators (see, e.g.,~\cite{babuvska1987hp}).
\end{itemize}
\noindent ~\qed
\newline \noindent 
While the energy norm cannot be computed exactly in  numerical verifications, it can be approximated 
by the error representation function. Hence, in Section~\ref{sec:verifications} when we present the 
error in the energy norm, it is the approximation of the energy norm through~\eqref{eq:energynorm_eqv}.

Since the KdV equation is a third order PDE, and we are interested in error bound in terms of $H^1$ and 
$L^2$, an application of the Aubin-Nitsche lift~\cite{aubin1987analyse,nitsche1972dirichlet} is required. To this end, we present the following lemma:
\begin{lem} \label{lem:sobolev_error}
Let $(u,q,p) \in \UUUT$ be the exact solution of the AVS-FE weak formulation~\eqref{eq:weak_form_KDV},  $(u^h,q^h,p^h) \in \UUUhT$ its corresponding AVS-FE approximation, and $m$ be the order of the 
Hilbert space in which we seek bounds. Then:
\begin{equation} \label{eq:sobolev_rate}
\ds \exists \, C > 0 \, : \norm{u-u^h}{H^m(\Omega_T)} \leq  C\, \ds h^{\,\nu},
\end{equation} 
where $h$ is the maximum element diameter, $\nu = \text{ min } (N - m,2N-S)$,  $N = \text{ min } (p_u +1,r)$, $p_u$ the minimum polynomial degree of approximation of $u^h$
in the mesh, $S$ the order of the PDE, i.e., $S=3$, and  $r$ the regularity of the solution $u$ of the governing PDE~\eqref{eq:KDV_eq}. 
\end{lem} 
\emph{Proof:}
The bound~\eqref{eq:sobolev_rate} is a consequence of:
\begin{itemize}
\item A quasi-best approximation property of the AVS-FE method in terms of the 
 norm on the trial space~\eqref{eq:broken_norms_KDV} (see Lemma 4.2.4 in~\cite{eirik2019thesis} for proof for a second order PDE);
\item The Aubin-Nitsche lift~\cite{aubin1987analyse,nitsche1972dirichlet,kastner2016isogeometric}
\item The existence of polynomial interpolation operators (see, e.g.,~\cite{babuvska1987hp}).
\end{itemize}
\noindent ~\qed
\newline \noindent 
The bound~\eqref{eq:sobolev_rate} depends not only on the approximation spaces and the FE mesh used 
but also on the order of the PDE. Thus, for polynomial degree of approximation $<2$, the $L^2$ and $H^1$ errors of $u-u^h$ are expected to converge at the same order.

\section{Numerical Verifications}  
\label{sec:verifications}

In this section, we present several numerical verifications of the AVS-FE method for the KdV equation. 
We consider both linear and nonlinear versions of the KdV equation and present academic problems to which exact solutions exist in an effort to confirm the rates of convergence predicted 
by the \emph{a priori} error estimates in Section~\ref{sec:avs-estimates}.
For the nonlinear problem we also present a verification of an $h-$adaptive algorithm.
The AVS-FE method is implemented using the saddle point interpretation~\eqref{eq:disc_mixed} into
the FE solver FEniCS~\cite{alnaes2015fenics}, which in turn employs the Portable, Extensible Toolkit for Scientific Computation 
(PETSc) library Scalable Nonlinear Equations Solvers (SNES)~\cite{abhyankar2018petsc,petsc-user-ref} to perform Newton iterations.

In all the following numerical verifications we use equal order continuous Lagrange polynomials for the 
trial variables $(u^h,q^h,p^h)$, whereas the error representation functions are discretized using 
discontinuous Lagrange polynomials of the same order.

\subsection{Linear KdV Equation}
\label{sec:lin_prob}
\begin{figure}[h!]
\subfigure[ \label{fig:line_p1_conv} Error norms. ]{\centering
 \includegraphics[width=0.45\textwidth]{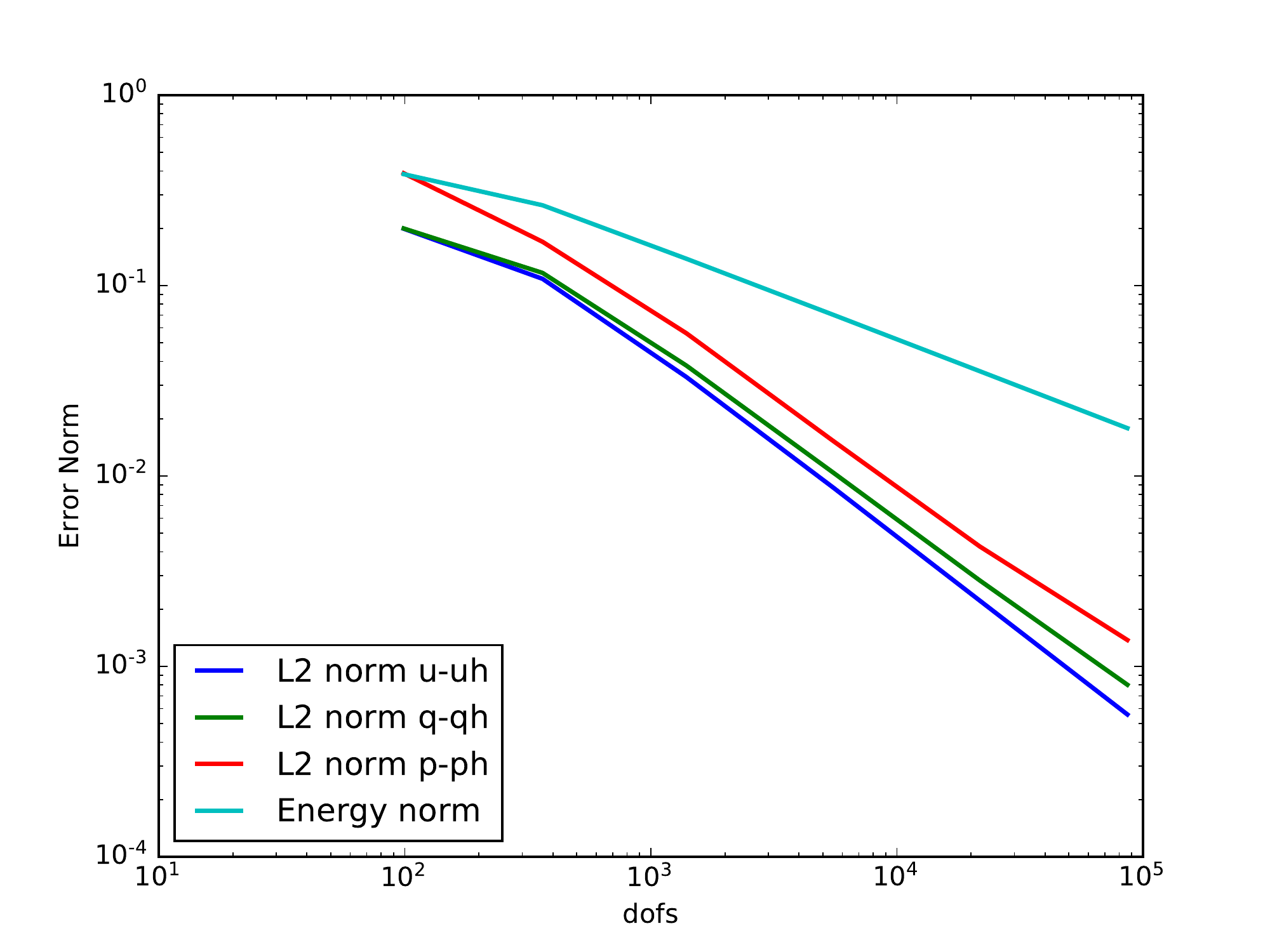}}
  \subfigure[ \label{fig:line_p1_rate} Rates of convergence.]{\centering
 \includegraphics[width=0.45\textwidth]{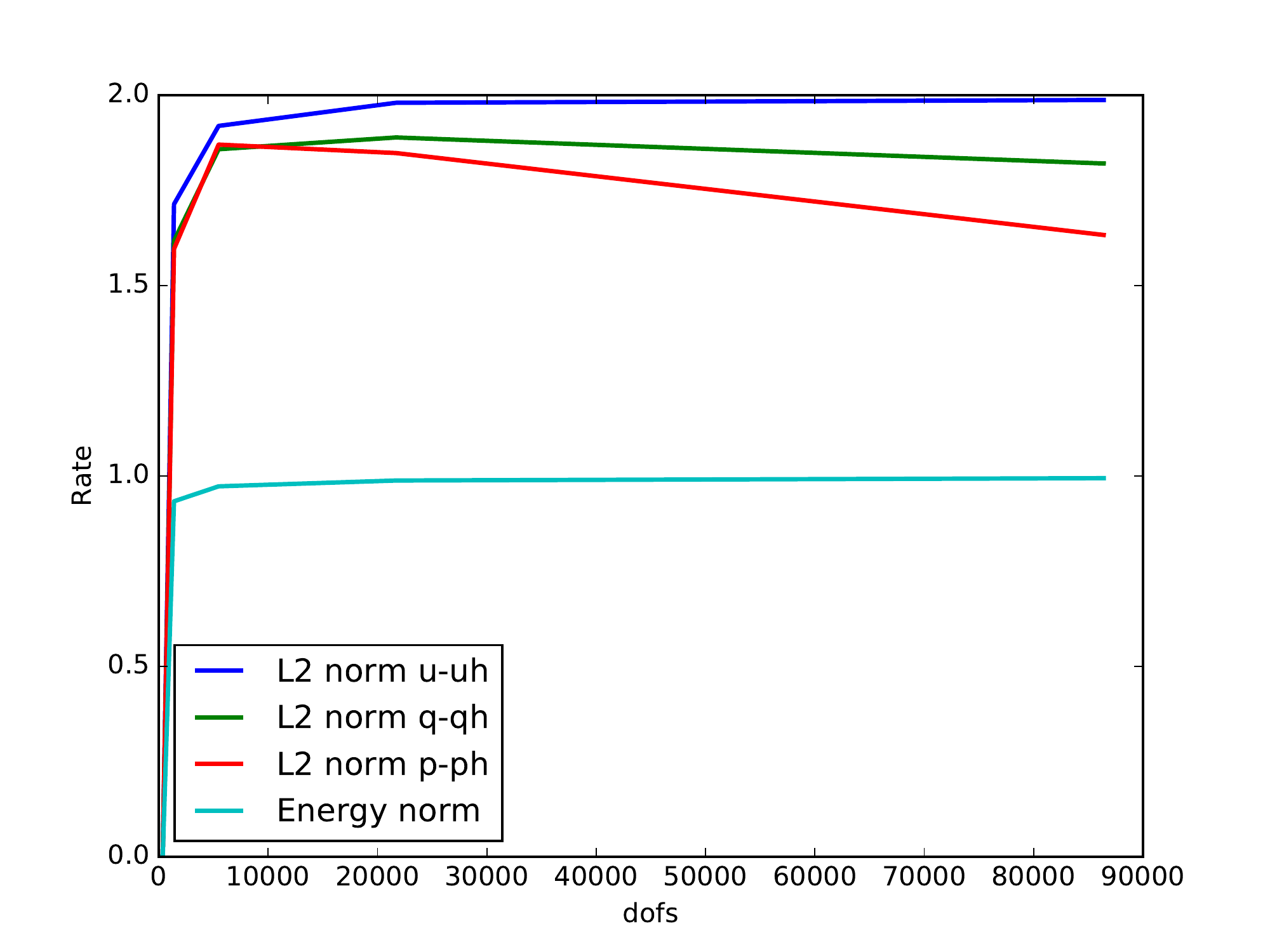}}
  \caption{\label{fig:Linear_p1}Error convergence results for uniform $h$-refinements for the linear KdV equation using linear polynomial approximations. }
\end{figure}

As an initial verification, we consider the linear version of~\eqref{eq:KDV_IBVP} introduced by Samii \emph{et al.} in~\cite{samii2016hybridized}, i.e., $\beta = 0$ and pick $\alpha = -1$,
the domain is $\Omega_T = ( 0,\pi ) \times (0,1.0s \rbrack$, and the 
 source is $f=0$.
As  initial and boundary conditions we use:
\begin{equation} \label{eq:lin_ex_BC}
\begin{array}{l}
  u(x,0)  =  \text{sin} \,x,  \\
  u(0,t)  =  -\text{sin} \,t,   \\  
  u(\pi,t)  =  \text{sin} \,t,   \\ 
  q(0,t)  =  \text{cos} \,t.  \\  
 \end{array}
\end{equation}
Thus, the exact solution is $\text{sin} (x-t) $.
First, we consider uniform mesh partitions of the space-time domain $\Omega_T$ that consists of triangular finite elements.  
In Figure~\ref{fig:line_p1_conv} we present the convergence plot for linear polynomial approximations for the $L^2$ norms of the approximation errors of all trial variables as well 
as the approximate energy norm. The corresponding rates of convergence are presented in Figure~\ref{fig:line_p1_rate} which reveals that the error $\norm{u-u^h}{\SLTOT}$ converges at a higher rate 
than predicted by~\eqref{eq:sobolev_rate}.

\begin{figure}[h!]
\subfigure[ \label{fig:line_p2_conv} Error norms. ]{\centering
 \includegraphics[width=0.45\textwidth]{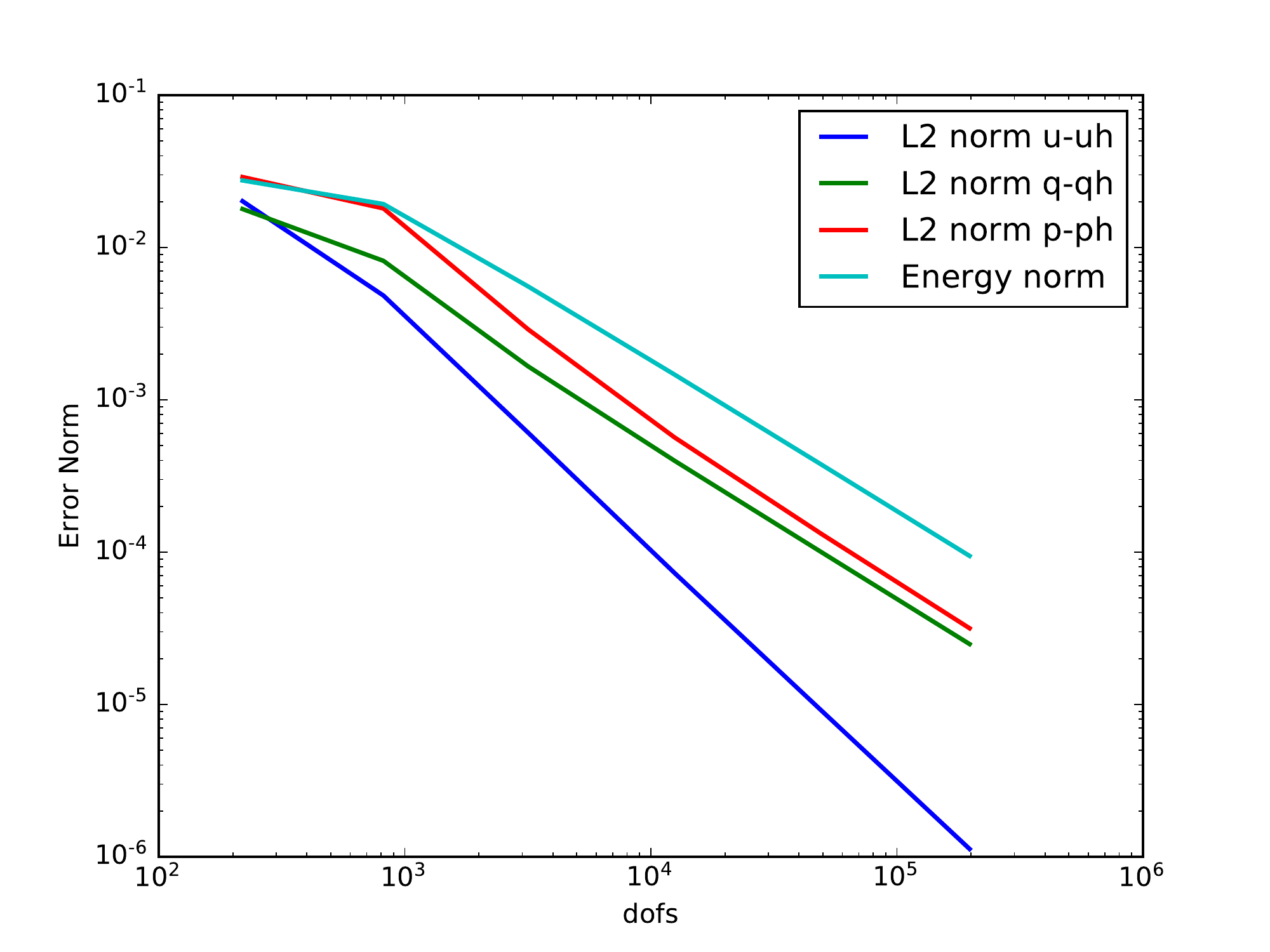}}
  \subfigure[ \label{fig:line_p2_rate} Rates of convergence.]{\centering
 \includegraphics[width=0.45\textwidth]{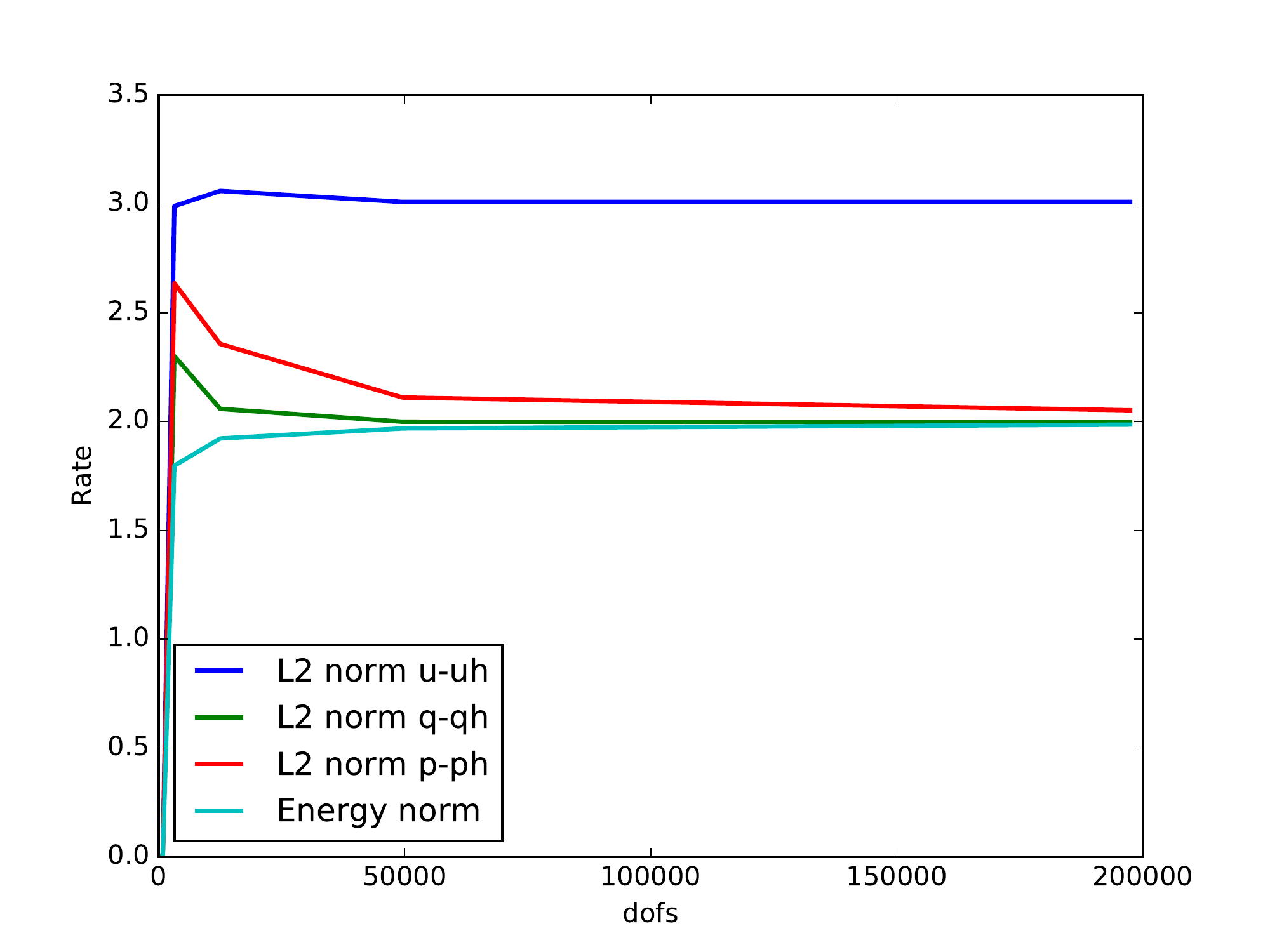}}
  \caption{\label{fig:Linear_p2}Error convergence results for uniform $h$-refinements for the linear KdV equation using quadratic polynomial approximations. }
\end{figure}

In Figure~\ref{fig:Linear_p2} we present similar results using quadratic polynomial approximations.
We do not present the convergence plots for $\norm{u-u^h}{\SHOOT}$ here but note that we have verified the rates predicted by~\eqref{eq:sobolev_rate} as well. 

\begin{figure}[h]
{\centering
\hspace{0.9in} \includegraphics[width=0.65\textwidth]{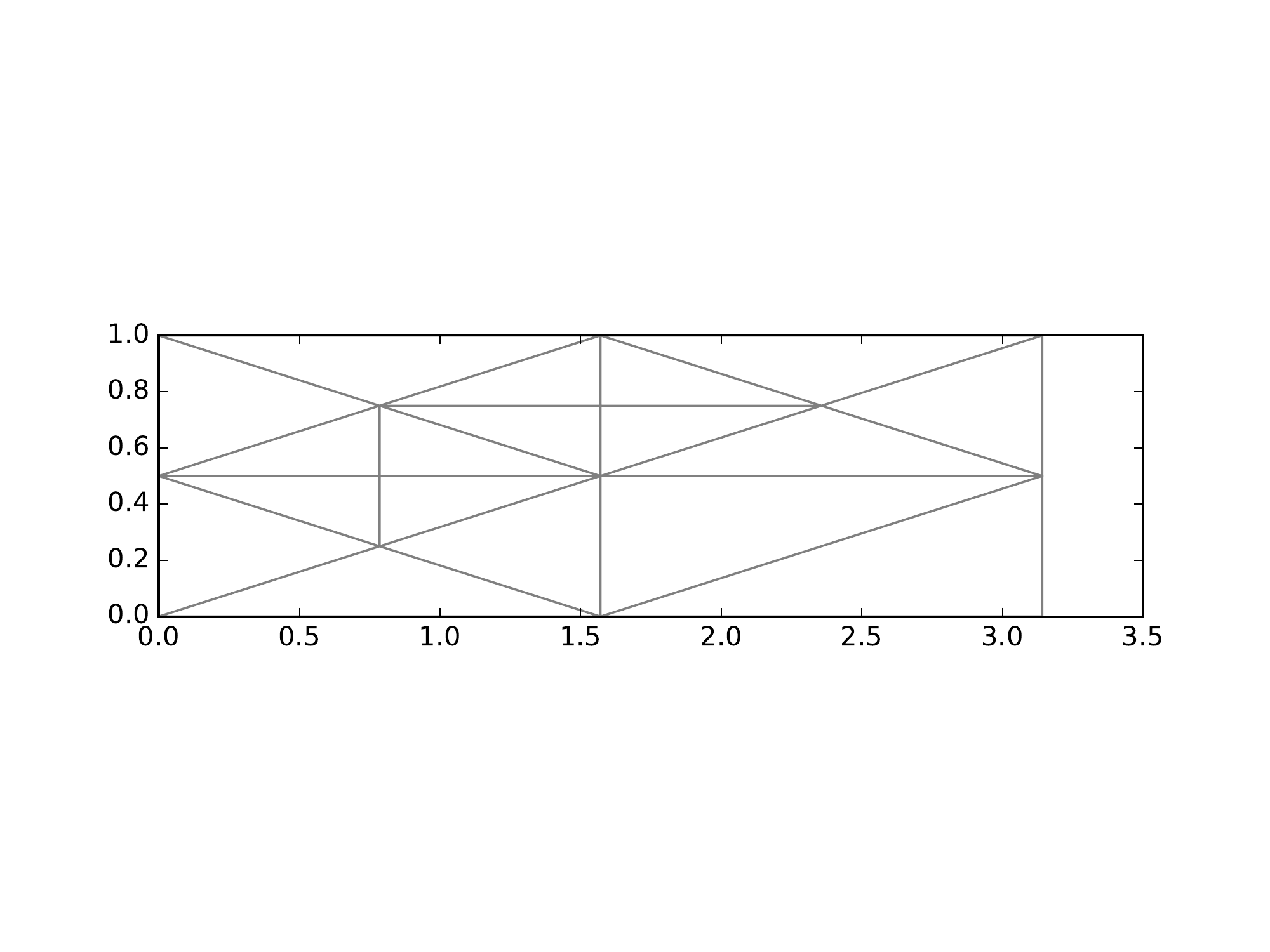}}
  \caption{\label{fig:unstruc_mesh} Example unstructured mesh.}
\end{figure}

As a final numerical verification for the linear KdV equation we consider non uniform meshes for the 
same problem data as the preceding example to ensure the convergence data is not influenced by the mesh 
structure. In Figure~\ref{fig:unstruc_mesh} an example of an
unstructured mesh used is shown.  The convergence properties for unstructured meshes remain unchanged 
as expected and we do not show further results for the linear case.

\subsection{Nonlinear KdV Equation}
\label{sec:nonllin_prob}
To study the convergence behavior for a nonlinear problem, we consider a case where  
$\beta = 3$, $\alpha = -1$, and domain $\Omega_T = ( 0,\pi )\times (0,0.5s \rbrack$.
The exact solution is chosen to be $u(x,t) = \text{sin} (\pi x)\,\text{sin} (\pi t)$, 
and we apply the differential operators of~\eqref{eq:KDV_eq} and~\eqref{eq:KDV_IBVP_first_order} to 
establish the source term $f(x,t)$ and boundary conditions on $u$ and $q$.
We first consider uniform meshes consisting of triangles to which we perform successive uniform 
mesh refinements.

\begin{figure}[h!]
\subfigure[ \label{fig:noline_p3_conv} Error norms. ]{\centering
 \includegraphics[width=0.45\textwidth]{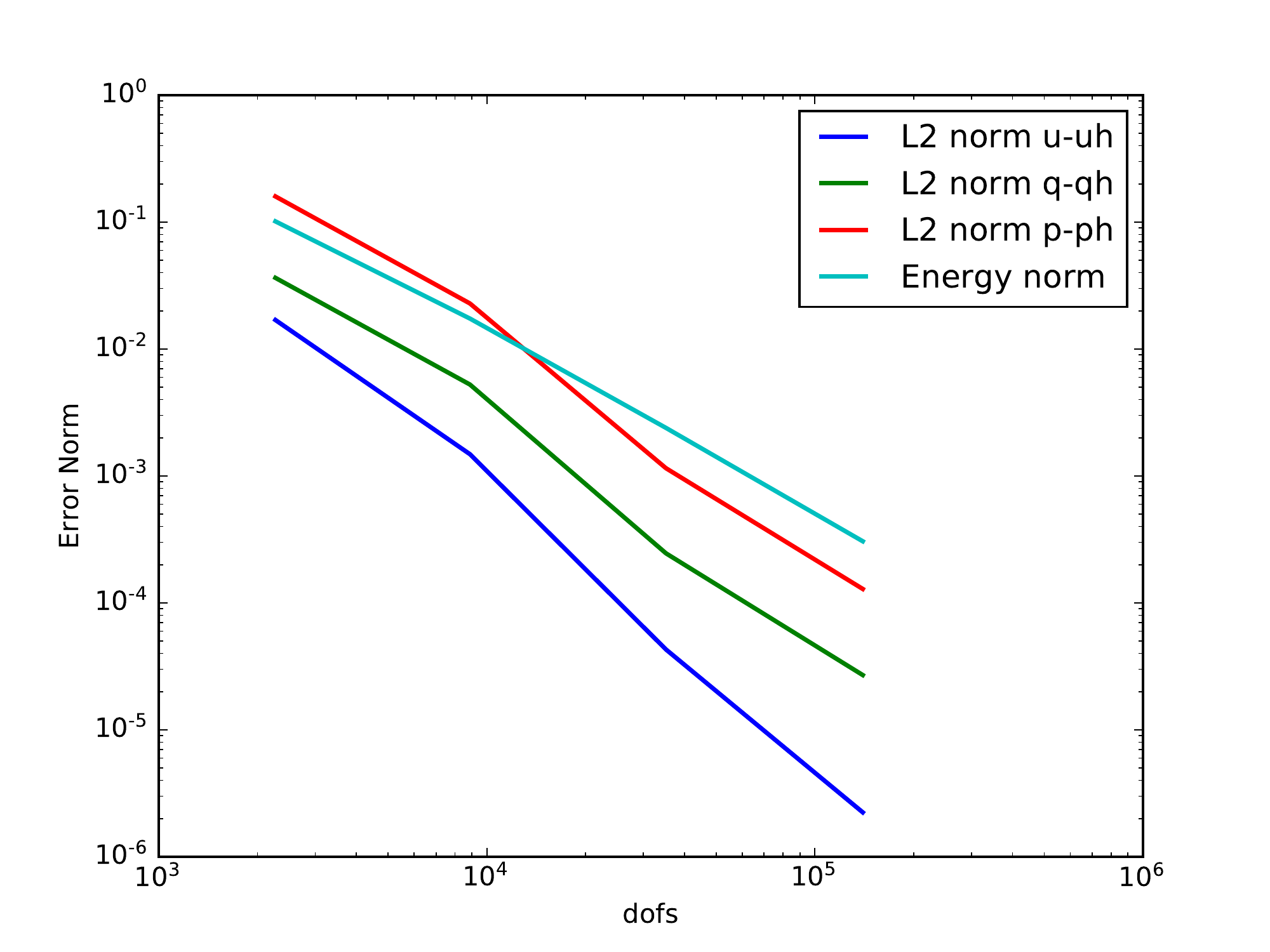}}
  \subfigure[ \label{fig:noline_p3_rate} Rates of convergence.]{\centering
 \includegraphics[width=0.45\textwidth]{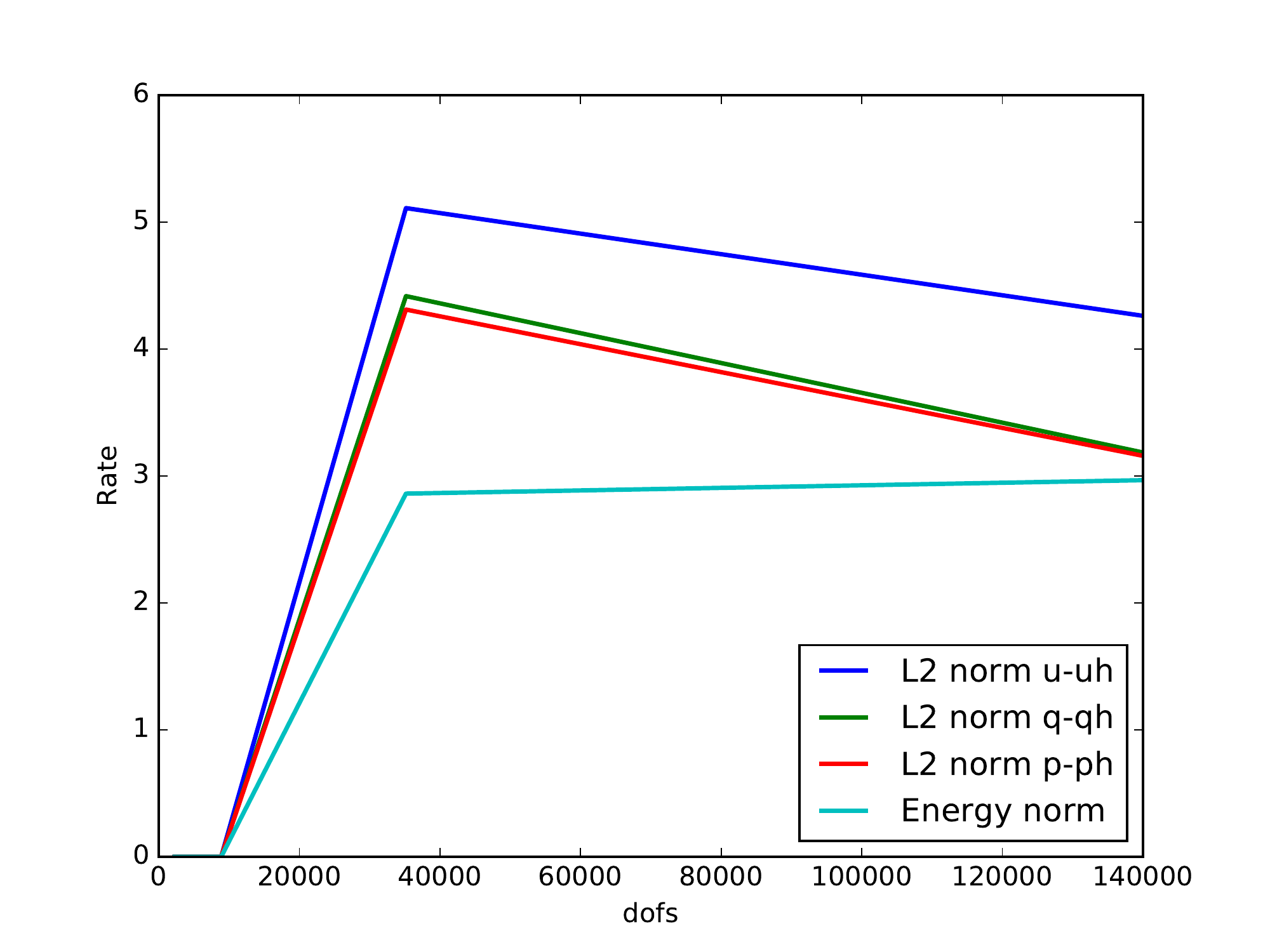}}
  \caption{\label{fig:NoLinear_p3}Error convergence results for uniform $h$-refinements for the nonlinear KdV equation using cubic polynomial approximations. }
\end{figure}

In Figure~\ref{fig:noline_p3_conv} we present the convergence plot for cubic polynomial approximations
 for the $L^2$ norms of the approximation errors of all trial variables as well 
as the approximate energy norm. The corresponding rates of convergence are presented in 
Figure~\ref{fig:noline_p3_rate} which reveals that the convergence rate of the error 
$\norm{u-u^h}{\SLTOT}$ approaches the rates predicted by~\eqref{eq:sobolev_rate} for the linear case.
Note that the rate of convergence of $\norm{u-u^h}{\SHOOT}$ agree with these predictions.
This is representative of a substantial number of numerical verifications for increasing degrees of approximation for uniform meshes.

\begin{figure}[h!]
\subfigure[ \label{fig:noline_p2_conv_ustruc} Error norms. ]{\centering
 \includegraphics[width=0.45\textwidth]{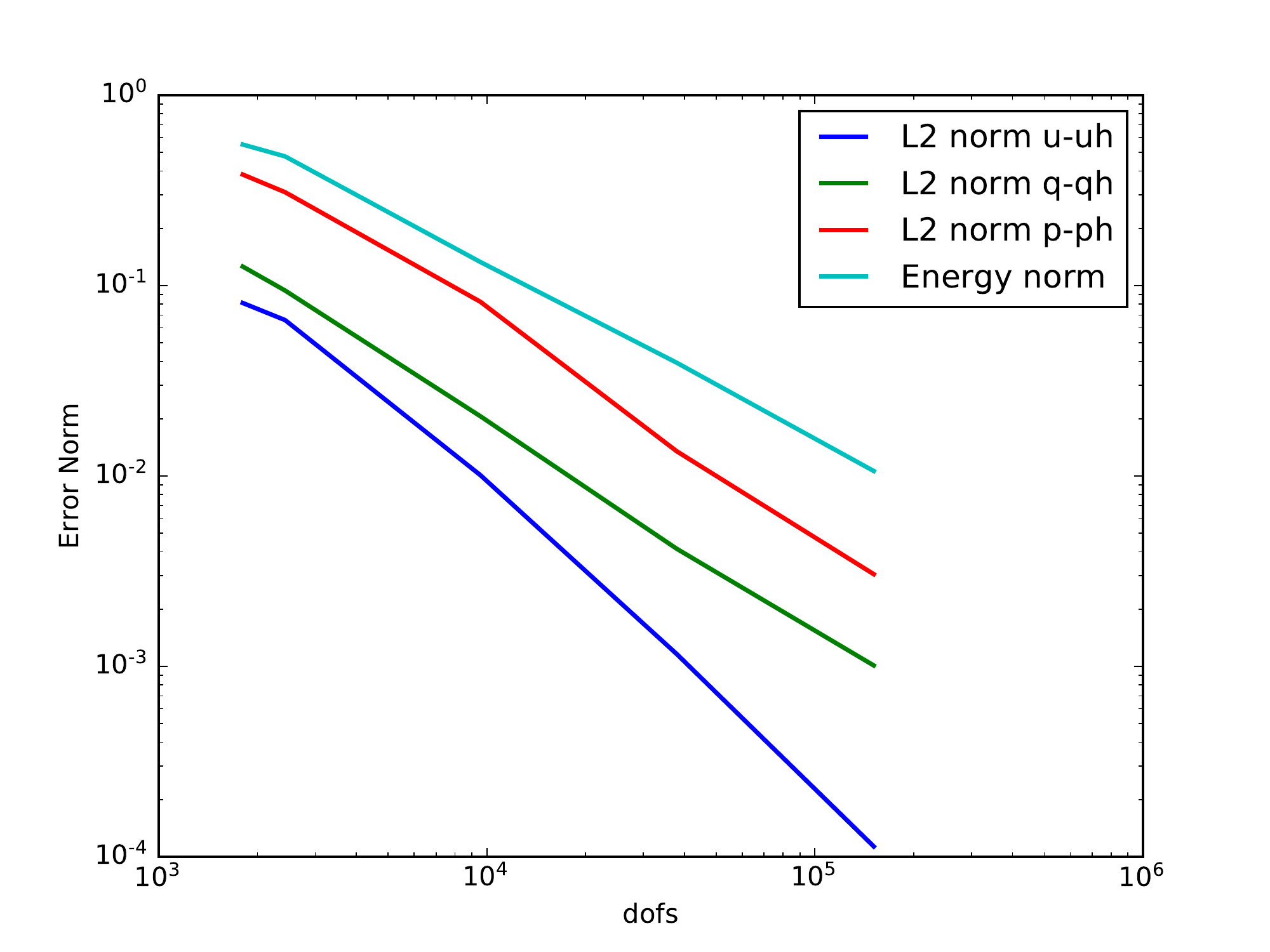}}
  \subfigure[ \label{fig:noline_p2_rate_ustruc} Rates of convergence.]{\centering
 \includegraphics[width=0.45\textwidth]{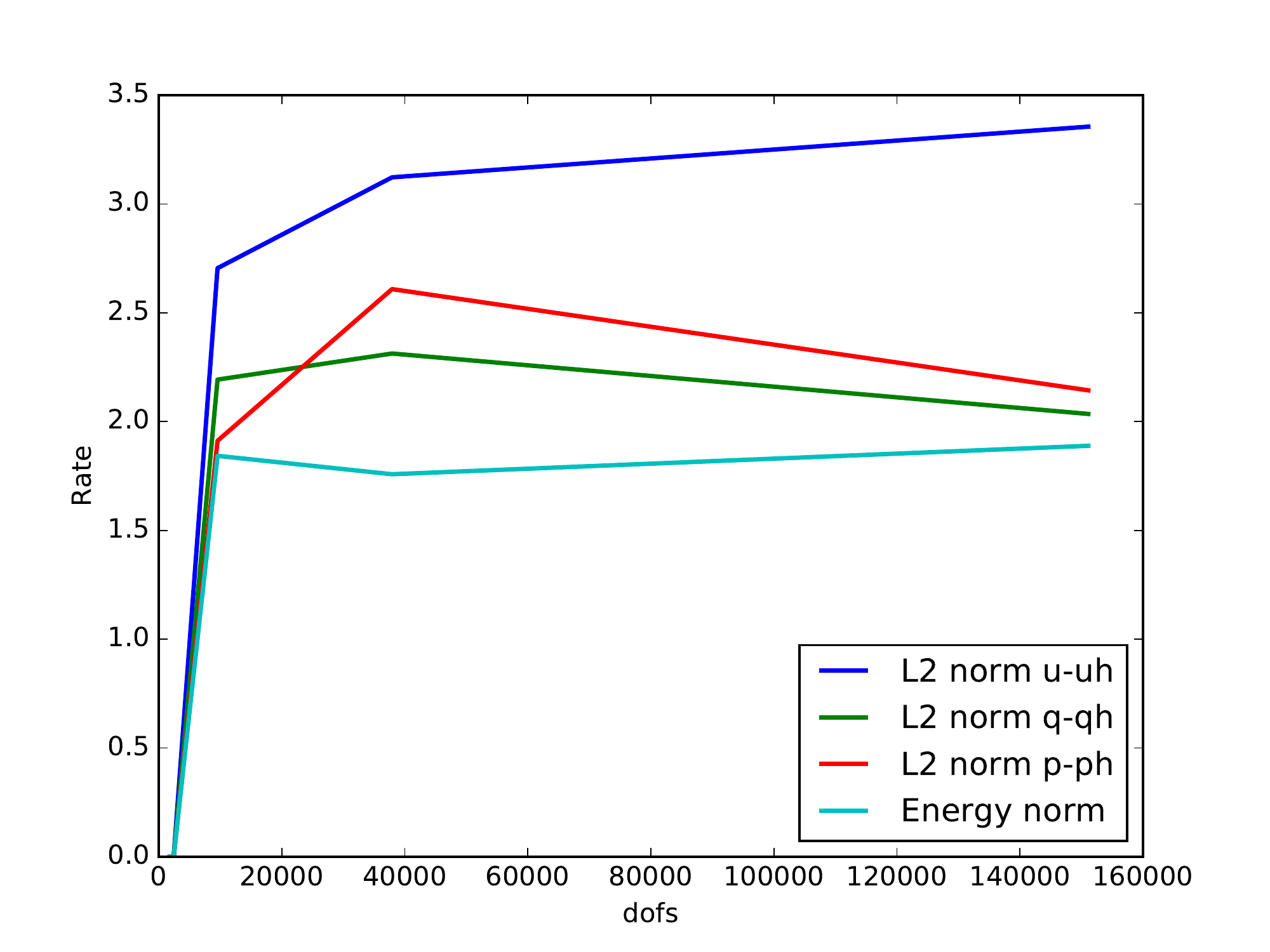}}
  \caption{\label{fig:NoLinear_p2_ustruc}Error convergence results for uniform $h$-refinements for the nonlinear KdV equation on a non uniform mesh using quadratic polynomial approximations. }
\end{figure}

As for the linear case, we again perform these verifications for non-uniform meshes. The convergence data from the case in which we use quadratic polynomial 
approximations is shown in Figure~\ref{fig:NoLinear_p2_ustruc}, where we see no effect on the convergence properties from the non uniformity of the meshes. 

\subsection{Adaptive Mesh Refinement}
\label{sec:adaptivity}

While the convergence behavior and optimal convergence rates provide confidence in the AVS-FE method applied 
to the KdV equation, uniform mesh refinements are generally not practical in physical applications.
Thankfully, the saddle point problem we solve~\eqref{eq:disc_mixed} also comes with "built-in" 
error estimators and indicators. We have used this estimator in the preceding verifications to 
estimate the energy norm. However, this norm can also be used to estimate the error \emph{a posteriori}
in cases where the exact analytic solution is unknown. Bounds on this estimator has been established 
for the DPG method in~\cite{Demkowicz2}, and its robustness has been numerically verified in numerous
papers. 
We therefore propose to employ the resulting error indicator - $\eta$, see~\eqref{eq:enrr_ind}, to drive mesh adaptive refinements according to
 the marking strategy and refinement criteria of D{\"o}rfler~\cite{dorfler1996convergent}.

\begin{figure}[h]
{\centering
\hspace{1in}\includegraphics[width=0.5\textwidth]{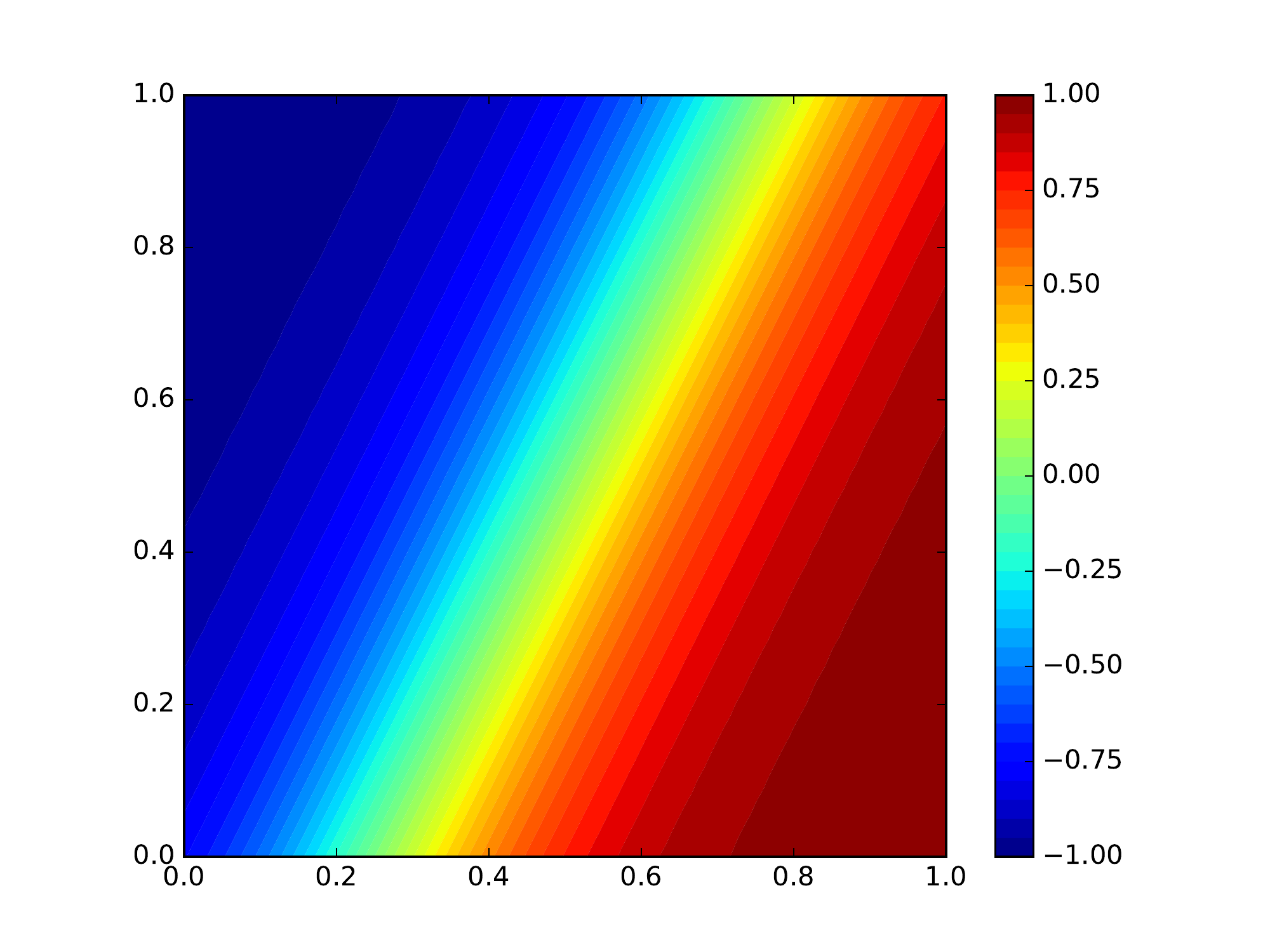}}
 \caption{\label{fig:exact} Exact solution $u(x,t)$ for the hyperbolic Tangent. }
\end{figure}

\begin{figure}[h!]
\subfigure[ \label{fig:in_mesh} Initial  uniform mesh. ]{\centering
 \includegraphics[width=0.45\textwidth]{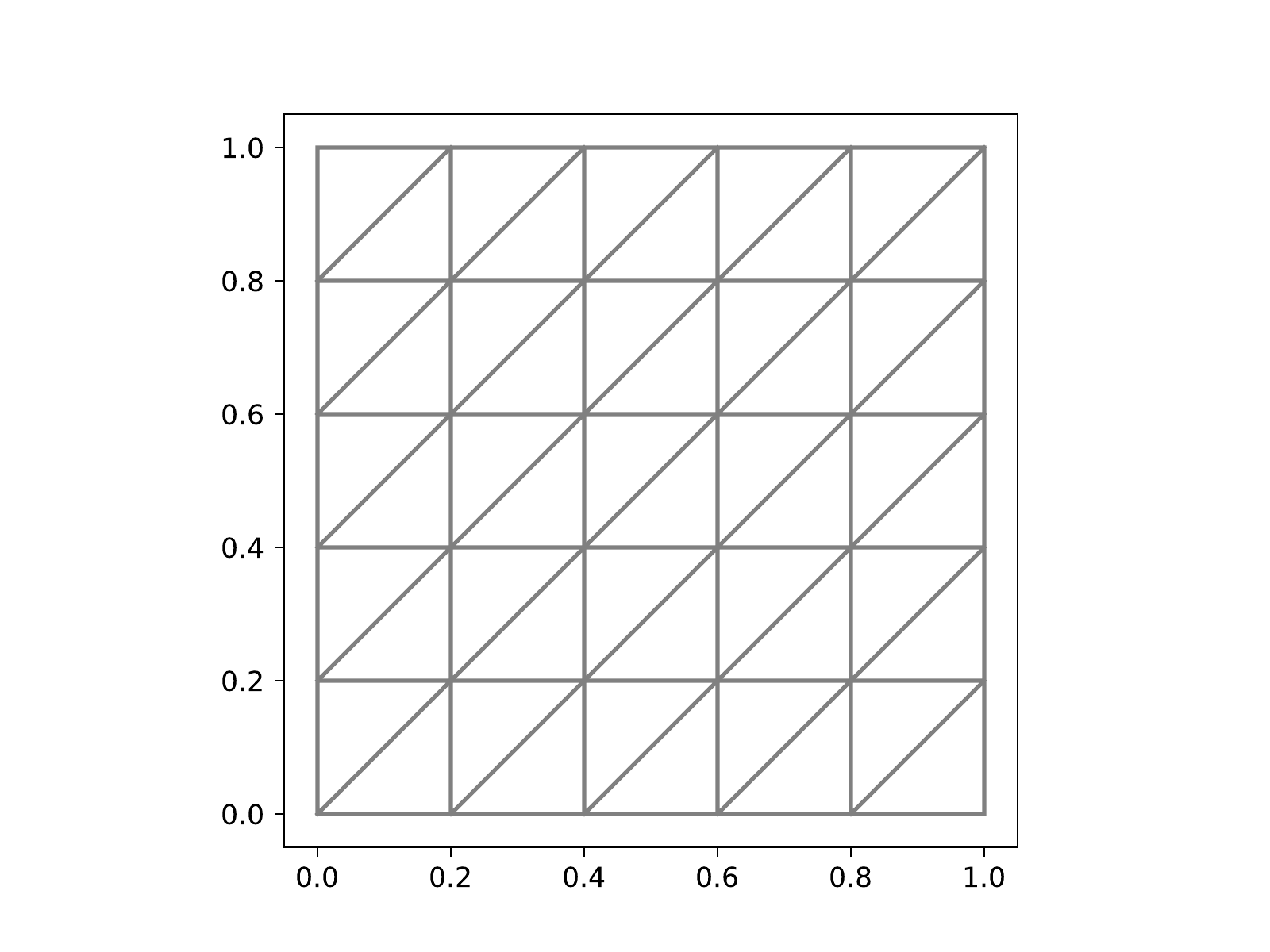}}
  \subfigure[ \label{fig:in_sol} Solution $u^h$.]{\centering
 \includegraphics[width=0.45\textwidth]{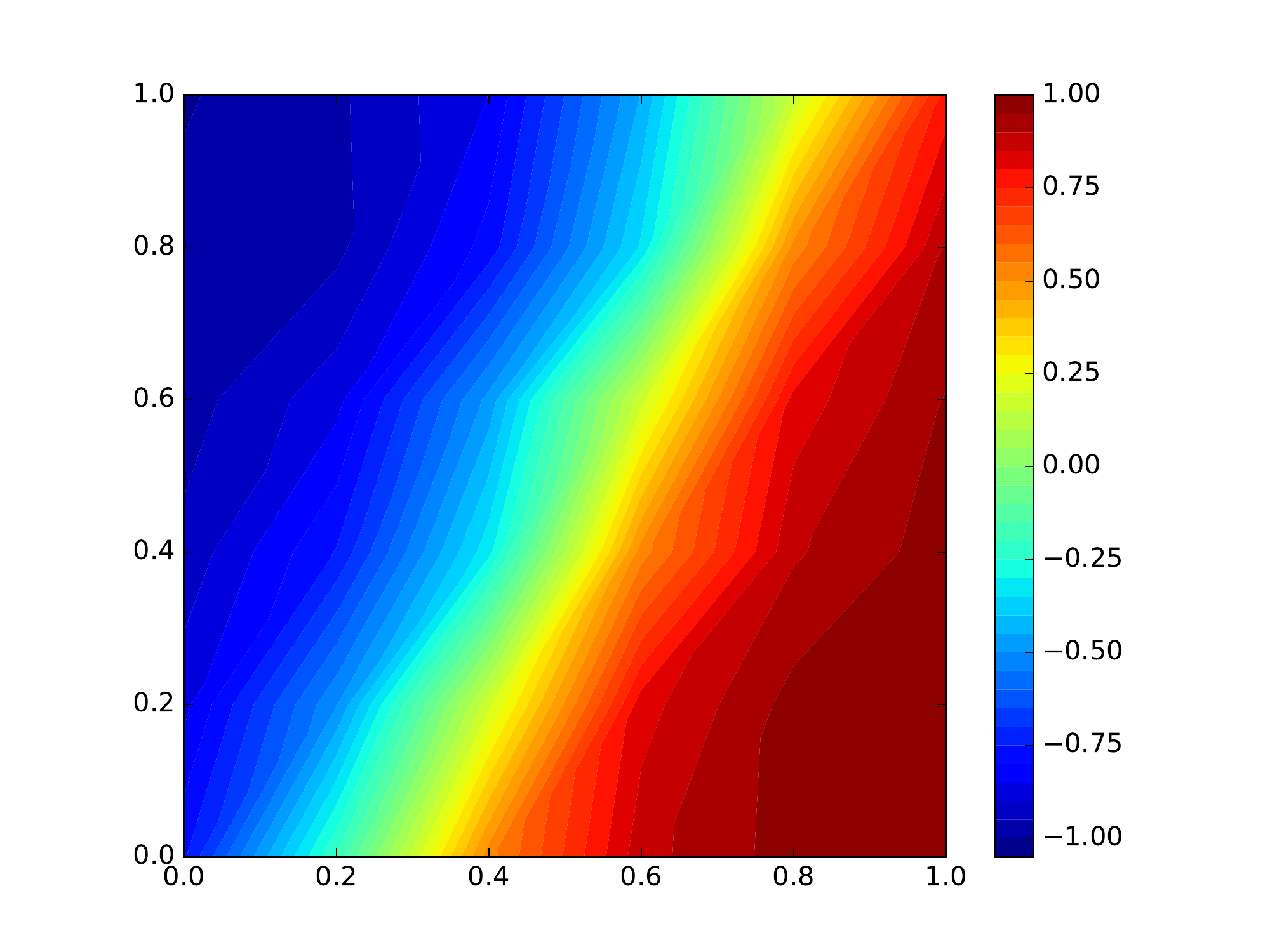}}
  \caption{\label{fig:initial_adaptive_step} Initial step.}
\end{figure}

\begin{figure}[h!]
\subfigure[ \label{fig:23_mesh} Adaptively refined mesh. ]{\centering
 \includegraphics[width=0.45\textwidth]{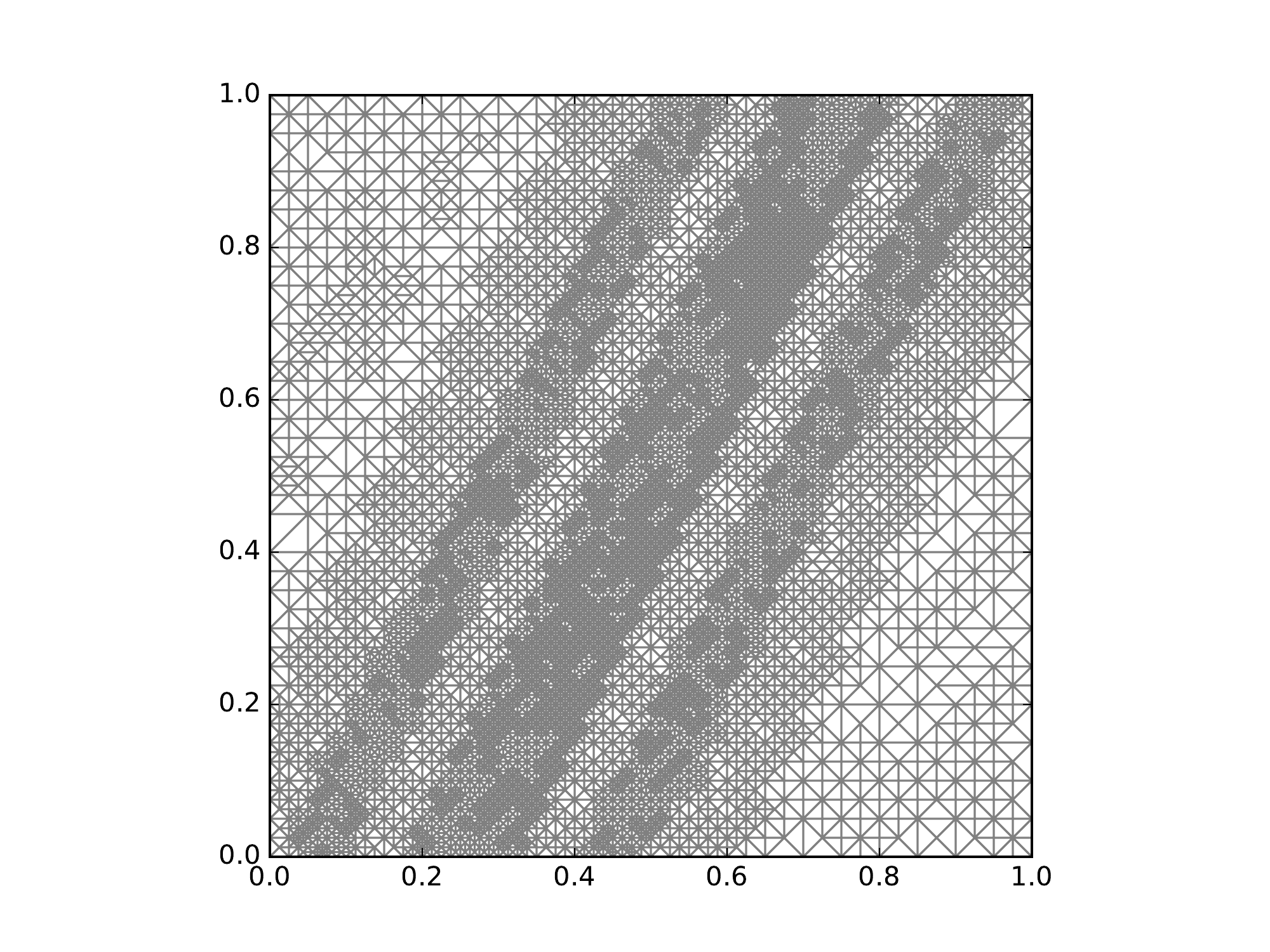}}
  \subfigure[ \label{fig:23_sol} Solution $u^h$.]{\centering
 \includegraphics[width=0.45\textwidth]{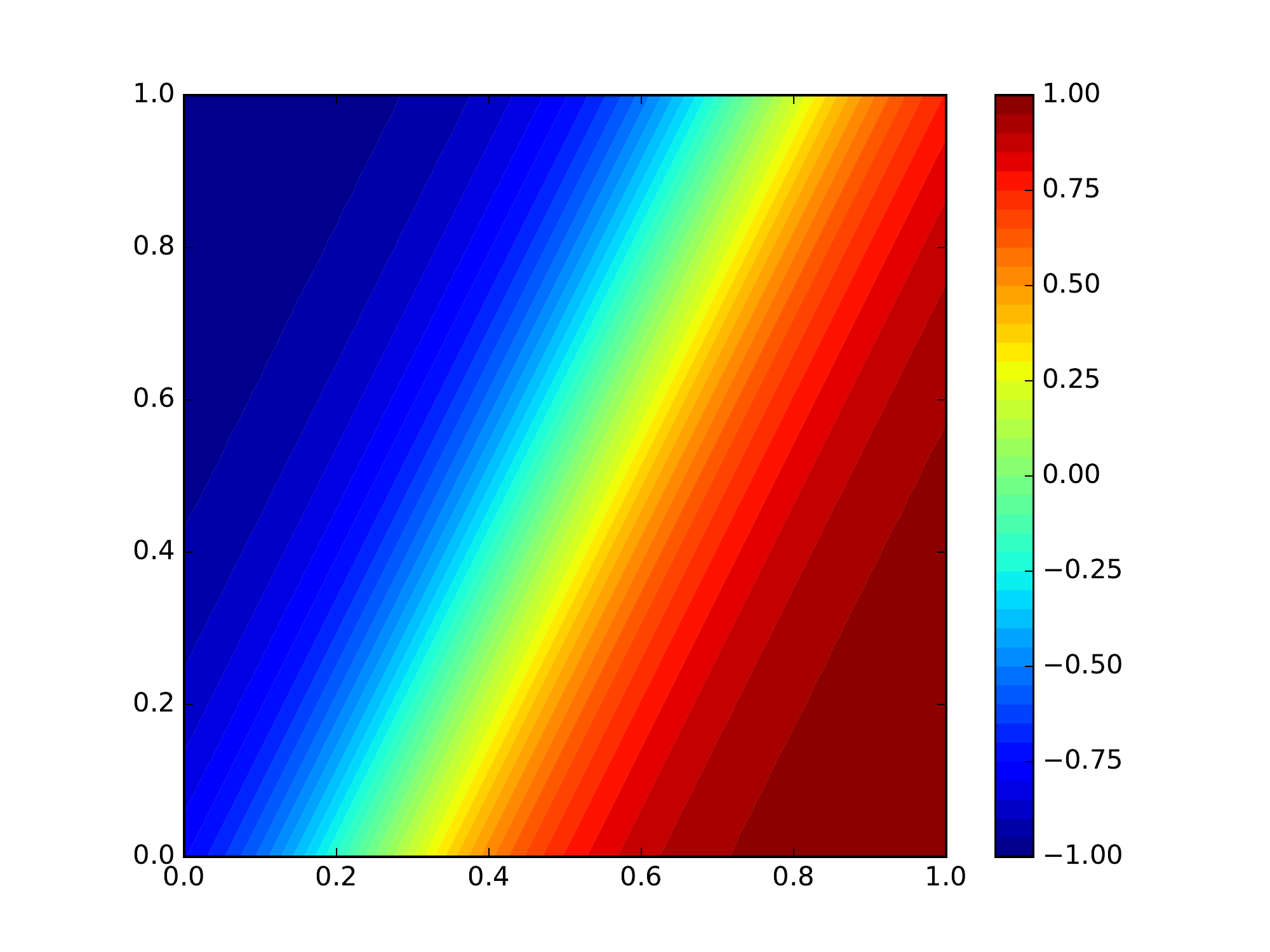}}
  \caption{\label{fig:23_adaptive_step} Final adaptive step.}
\end{figure}

\begin{figure}[h]
{\centering
\hspace{0.9in} \includegraphics[width=0.65\textwidth]{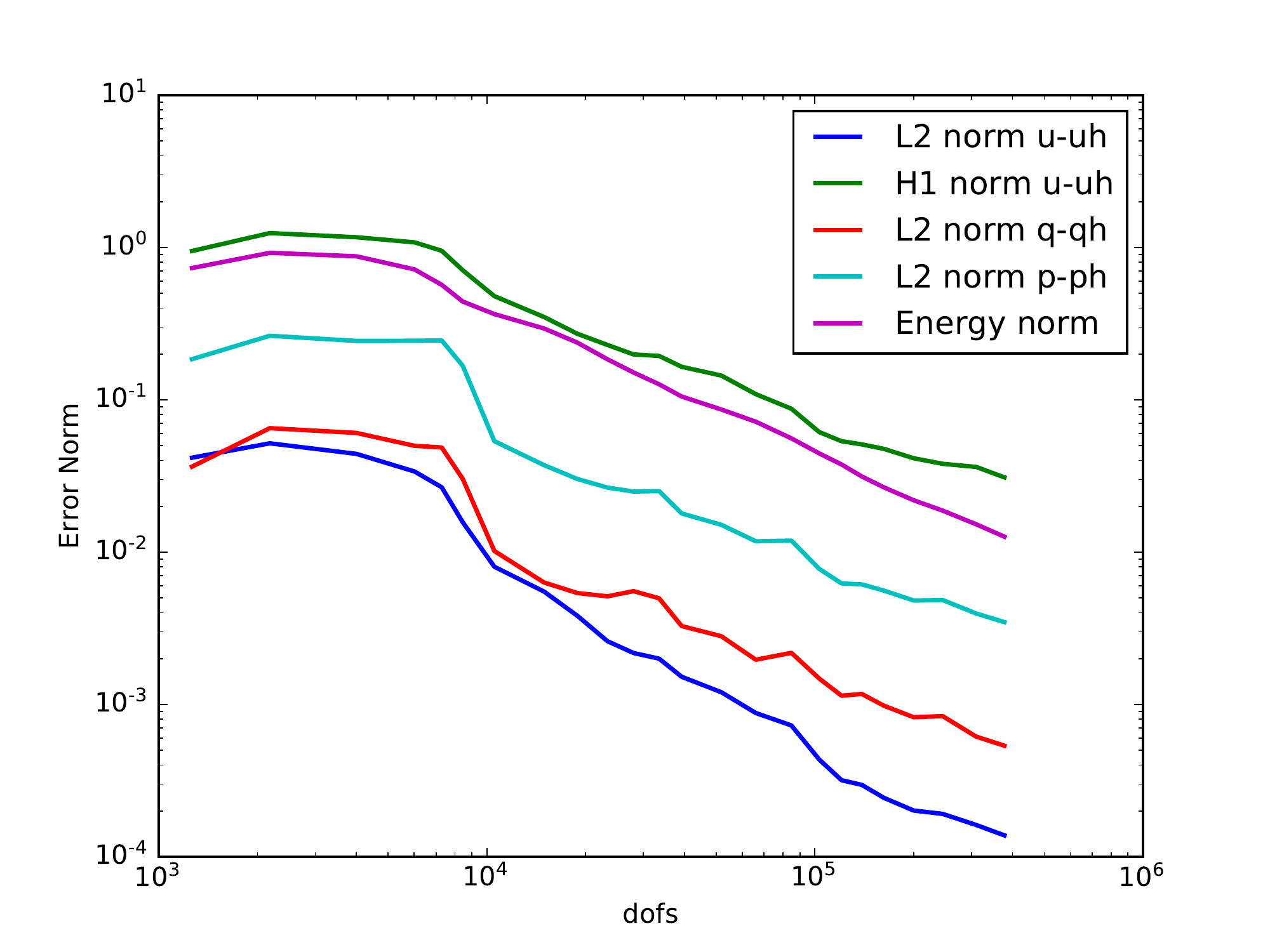}}
  \caption{\label{fig:convergence_results_adaptive} Error convergence results for adaptive $h$-refinements.}
\end{figure}

We still consider the nonlinear KdV equation with
$\beta = 3$, $\alpha = -1$, the domain is now $\Omega_T = (0,1 ) \times (0,1.0s \rbrack$,
and we use quadratic polynomial approximations.
We pick the exact solution to be the hyperbolic Tangent function:
\begin{equation}
\label{eq:prop_fron}
\begin{array}{c}
\ds u(x,t) = \text{tanh} \left( \frac{\ds x - 0.5t -0.25}{\ds 0.25  } \right),
\end{array}
\end{equation}
from which we establish the source $f(x,t)$ and boundary conditions on $u$ and $q$.
This exact solution is shown in Figure~\ref{fig:exact}.
In Figures~\ref{fig:initial_adaptive_step} 
and~\ref{fig:23_adaptive_step} we show the first and final meshes and 
corresponding solutions from the refinement process. 
Comparison of the first mesh in Figure~\ref{fig:in_mesh} and the final mesh in Figure~\ref{fig:23_mesh}
clearly show that the built in error indicator performs very well. 
The mesh refinements are focused along the diagonal where the changes in the 
hyperbolic Tangent are the greatest (see Figure~\ref{fig:exact}). 
At the same time, the elements far away from this region are far less refined further indicating 
the accuracy of the error indicators.
Lastly, we show the convergence history of the adaptive refinement process in
Figure~\ref{fig:convergence_results_adaptive}.

\section{Conclusions}
\label{sec:conclusions}

We have introduced an unconditionally stable FE method for the KdV equation,
the AVS-FE method. 
Highly accurate FE approximations are established in both space and time due to the unconditional stability of the AVS-FE method. 
We use a saddle point representation of the method which allows us to 
implement it into the high-level FE solver FEniCS~\cite{alnaes2015fenics} in 
a straightforward manner.

We have presented  \emph{a priori} error bounds for the linear version of the KdV which we use 
to predict convergence rates. 
The AVS-FE approximations of the linear version of the KdV follow the predictions for quadratic 
polynomial approximations and higher. Whereas for linear approximations, the rate of convergence is higher than expected, and is shown in Section~\ref{sec:lin_prob}. 
For the nonlinear version, the AVS-FE approximations converge at the expected rates for  
polynomial approximations and is presented in Section~\ref{sec:nonllin_prob}. For both linear and nonlinear versions, the convergence behavior is unaffected by unstructured mesh partitions. 
This convergence behavior is noteworthy, as we are able to achieve optimal rates of convergence throughout the space-time domain for any degree of approximation. This is in contrast to classical 
time stepping techniques which require investigation of individual time steps to ascertain the
convergence behavior.
In Section~\ref{sec:lin_prob}, we solve a problem introduced by Samii \emph{et al.}~\cite{samii2016hybridized}. 
A direct comparison between the results for the linear KdV equation in Section~\ref{sec:lin_prob}  and those in~\cite{samii2016hybridized} is not straightforward. However, we note that the rates of convergence achieved for $\norm{u-u^h}{\SLTO}$ are the same at the reported times in~\cite{samii2016hybridized}.

While the space-time convergence properties gives confidence in our numerical scheme, a major   strength of the AVS-FE method is the built-in error estimator and indicator. 
Hence, we are able to establish an adaptive numerical scheme at no additional computational cost making 
the AVS-FE highly competitive with classical methods. 
Furthermore, the unconditional stability property of the AVS-FE allows us to compute FE 
approximations on very coarse initial meshes and then employ mesh adaptive processes 
to establish highly accurate FE approximations.
As shown in Section~\ref{sec:adaptivity}, this allows us to implement mesh adaptive algorithms 
which are capable of driving the approximation error towards zero. 
In future efforts we also plan to employ alternative error estimates in terms of local quantities of 
interest as proposed in~\cite{valseth2020goal}.

The AVS-FE method is capable of delivering highly accurate
space-time computations for surface waves in water as shown in this presentation, and we plan to
investigate more complex physical phenomena in future works.

\section*{Acknowledgements}
This work has been supported by the NSF PREEVENTS Track 2 Program,
under  NSF Grant Number  1855047.

%The authors are grateful for the contributions of ...

\newpage
%=======================
% Bibliography
%=======================]
\bibliographystyle{elsarticle-num}
 \bibliography{references_eirik}
\end{document}